\documentclass[11pt]{amsart}
\usepackage{graphicx}
\usepackage{color}
\usepackage{latexsym}
\usepackage{amscd}
\usepackage{a4wide}
\usepackage{epsfig}
\usepackage{float}
\usepackage{amsmath, amssymb, amsthm, amsfonts, amsgen} 
\usepackage[outdir=./]{epstopdf}

\graphicspath{ {./} }

\usepackage{esint}

\usepackage{enumitem}
\usepackage{cite}

\usepackage{color}

\renewcommand{\[}{\begin{equation*}}
\renewcommand{\]}{\end{equation*}}
\newcommand{\grad}{\ensuremath{\nabla}}
\newcommand{\p}{\ensuremath{\partial}}
\newcommand{\ve}{\ensuremath{\varepsilon}}
\newcommand{\vp}{\ensuremath{\varphi}}
\newcommand{\al}{\ensuremath{\alpha}}
\newcommand{\be}{\ensuremath{\beta}}
\newcommand{\ga}{\ensuremath{\gamma}}
\newcommand{\ta}{\ensuremath{\tau}}
\newcommand{\ka}{\ensuremath{\kappa}}

\renewcommand{\ss}{\ensuremath{\subset}}
\newcommand{\sm}{\ensuremath{\setminus}}
\newcommand{\ra}{\ensuremath{\rightarrow}}
\newcommand{\mt}{\ensuremath{\mapsto}}

\newcommand{\mb}[1]{\ensuremath{\;\mbox{#1}\;}}
\newcommand{\cd}{\ensuremath{\cdot}}
\newcommand{\ot}{\ensuremath{\otimes}}
\newcommand{\ti}{\ensuremath{\times}}

\newcommand{\xb}{\ensuremath{\bar{x}}}

\newcommand{\wt}[1]{\widetilde{#1}}
\newcommand{\ol}[1]{\overline{#1}}

\newcommand{\cof}{\textnormal{cof}\;}
\renewcommand{\div}{\textnormal{div}\;}
\newcommand{\adj}{\textnormal{adj}}

\newcommand{\dist}{\textnormal{dist}\;}
\newcommand{\supp}{\textnormal{supp}\;}

\newcommand{\R}{\mathbb{R}}

\newcommand{\N}{\mathbb{N}}

\newcommand{\bb}[1]{\ensuremath{{\mathbb{#1}}}}

\newcommand{\D}{\Delta}
\renewcommand{\d}{\delta}
\newcommand{\A}{\mathcal{A}}

\newcommand{\La}{\mathcal{L}}

\newcommand{\Om}{\Omega}

\newtheorem{de}{Definition}[section]
\newtheorem{thm}[de]{Theorem}

\newtheorem{lem}[de]{Lemma}

\newtheorem{re}[de]{Remark}
\allowdisplaybreaks
\begin{document}
\title{Everywhere regularity results for a polyconvex functional in finite elasticity} 
\author[Marcel Dengler]{M. Dengler}
\address[M. Dengler]{Fliederweg 1, 72189 Vöhringen, Germany.}
\email{marci.dengler@web.de}

\maketitle 
\begin{abstract}
\noindent Here we develop a regularity theory for a polyconvex functional in $2\ti2-$dimen-sional compressible finite elasticity. In particular, we consider energy minimizers/stationary points of the functional
\begin{equation*}I(u)=\int\limits_\Om{\frac{1}{2}|\grad u|^2+\rho(\det\grad u)\;dx},\label{eq:SA.1.1}\end{equation*}
where $\Omega\ss\R^2$ is open and bounded, $u\in W^{1,2}(\Om,\R^2)$ and $\rho:\R\ra\R_0^+$ smooth and convex with $\rho(s)=0$ for all $s\le0$ and $\rho$ becomes affine when $s$ exceeds some value $s_0>0.$ Additionally, we may impose boundary conditions.\\

\noindent The first result we show is that every stationary point needs to be locally Hölder-continuous. 
Secondly, we prove that if $\|\rho'\|_{L^\infty(\R)}<1$ s.t.\! the integrand is still uniformly convex, then all stationary points have to be in $W_{loc}^{2,2}.$ Next, a higher-order regularity result is shown. 
Indeed, we show that all stationary points that are additionally of class $W_{loc}^{2,2}$ and whose Jacobian is suitably Hölder-continuous are of class $C_{loc}^{\infty}.$ As a consequence, these results show that in the case when $\|\rho'\|_{L^\infty(\R)}<1$ all stationary points have to be smooth.  \\

\end{abstract}

\section{Introduction}
For $\Om\ss \R^2$ open, bounded, and with $\p\Om\in C^{0,1}$ define the functional $I:W^{1,2}(\Om,\R^2)\rightarrow \R$ by
\begin{equation}
I(u):=\int\limits_{\Om}{\frac{1}{2}|\grad u|^2+\rho(\det \grad u)\; dx}
\label{eq:1.1}
\end{equation}
for all $u\in W^{1,2}(\Om,\R^2).$ The function $\rho\in C^{\infty}(\R)$ is defined by
\begin{equation}
\rho(s)=\left\{\begin{array}{ccc}
0& {\mbox{if}}& s\le0,\\
\rho_1(s)& {\mbox{if}}& 0\le s\le s_0,\\
\gamma s+\ka&{\mbox{if}}& s_0\le s,
\end{array}
\right.
\label{eq:1.2}
\end{equation}
for some constants $\gamma>0,$ $s_0\ge0$ and $\ka\ge-\ga s_0.$ 
Here $\rho_1:[0,s_0]\rightarrow \R$ is a smooth and convex function on $[0,s_0]$ satisfying the boundary conditions $\rho_1(0)=0$ and $\rho_1(s_0)=\gamma s_0+\ka$ and the connections need to be in such a way that $\rho$ is smooth everywhere.\footnote{It would be enough to assume $\rho\in C^k$ for some $k\ge2$ with the necessary changes in the statements that have to be made. The results below also remains true for variations of the integrand for example one can pick $\rho$ to be of the form
\begin{equation*}
\rho(s)=\left\{\begin{array}{ccc}
\gamma |s|+\ka&{\mbox{if}}& s\le -s_0,\\
\rho_1(s)& {\mbox{if}}& -s_0\le s\le s_0,\\
\gamma s+\ka&{\mbox{if}}& s_0\le s,
\end{array}\right.
\end{equation*} 
where $\rho_1:[-s_0,s_0]\rightarrow \R_0^+$ is again a smooth and convex function on $[-s_0,s_0]$ satisfying the conditions $\rho_1(0)=0,$ $\rho_1(s_0)=\gamma s_0+\ka$ and $\rho_1(-s_0)=\gamma|s_0|+\ka$ with smooth connections.}  Note that $\rho$ is convex on the whole real line. Hence, the complete integrand is polyconvex. Recall, again, we call $f:\R^{2\ti2}\ra\ol{\R}$ polyconvex, if there exists a convex function $g:\R^5\ra\ol{\R}$ s.t. $f(\xi)=g(\xi,\det\xi)$ for all $\xi\in\R^{2\ti2}.$
The behaviour of the functional depends mainly on the parameter $\ga.$ If $\ga\ra0$ then the functional turns into the well known Dirichlet energy. In the regime $0<\ga<1$ the functional is uniformly convex. If $\ga\ge1$ then the functional is genuinely polyconvex.
Recall that for $u\in W^{1,2}(\Om,\R^2)$ the Jacobian is in $L^1,$ i.e. $\det\grad u\in L^1(\Om)$ so $\rho(\det\grad u(\cd))\in L^1(\Om)$ and hence $I(u)$ is well defined for all $u\in W^{1,2}(\Om,\R^2).$\\

Furthermore, we introduce the set of admissible functions
\[\A_{u_0}:=\{u\in W^{1,2}(\Om,\R^2): u=u_0 \;\mb{on}\; \p\Om\},\]
where the boundary condition must be understood in the trace sense.\\

Then the functional $I$ attains its minimum in $A_{u_0}.$ A proof is classically obtained via the direct method of the Calculus of Variations and \cite[Prop.A.3] {BM84}.\\

\textbf{Notation:} Let $A,B\in \R^{m_1\times\cdots\times m_n},$ where $m_k\in\N\sm\{0\}$ for all $0\le k\le n.$ We will denote the Frobenius inner product by $A\cdot B:=\sum\limits_{i_1,\ldots, i_n}A_{i_1,\ldots, i_n}B_{i_1,\ldots, i_n}$ and the corresponding norm by $|A|:=(A\cdot A)^{\frac{1}{2}}$ for all $A,B\in \R^{m_1\times\cdots\times m_n}.$
For a matrix $A\in \R^{n\ti n}$ we denote the determinant by $d_A:=\det A$ (when there is no confusion, we will suppress the dependence on $A$.) and the cofactor by $\cof A:=(\adj \;A)^T.$ It is well known, that the cofactor is the derivative of the determinant, i.e. $\grad_A(\det A)=\cof A.$ Moreover, for a $2\times2-$matrix $A$ the cofactor takes the simple form 
\[
\cof A=\begin{pmatrix}a_{22} &-a_{21}\\-a_{12} & a_{11}\end{pmatrix}
\]
which is a rearrangement of the original matrix (up to signs). Note, that $|\cof A|=|A|.$\\
For two vectors $a\in \R^n,b\in \R^m$ we define the tensor product $a\ot b\in \R^{n\ti m}$ by $(a\ot b)_{i,j}:=(ab^T)_{i,j}=a_ib_j$ for all $1\le i\le n,$ $1\le j\le m.$ $B(x,r):=\{y\in\bb{R}^2:|y-x|<r\}$ be the open ball with center $x$ in $\R^2.$ If there is no confusion, we use the abbreviation $B_r.$ We want $\mathcal{H}^s$ to be the $s-$dimensional Hausdorff measure for any $s\in \R_0^+,$ and $\La^2$ to be the $2-$dimensional Lebesgue measure, as usual, we use $dx=d\La^2$ for short.\\

\begin{de}[stationary point]
We call $u\in\A_{u_0}$ a stationary point of the functional \eqref{eq:1.1} if it satisfies the Euler-Lagrange equation (ELE) in the weak form, which is given by
\begin{equation}
\int\limits_{\Om}{(\grad u+\rho'(d_{\grad u})\cof \grad u) \cdot\grad\vp\; dx}=0\; {\mbox{for all}}\; \vp\in W_0^{1,2}(\Om,\R^2).
\label{eq:ELE1.1}
\end{equation}
\end{de}
\vspace{0.3cm}

Firstly, we show that any stationary point of the functional \eqref{eq:1.1} needs to be locally Hölder continuous.
\begin{thm}[Hölder continuity]
Let $0< \ga<\infty$ and suppose that $u\in\A_{u_0}$ be a stationary point of the functional $I$ as defined in \eqref{eq:1.1}. Then $u$ is locally Hölder continuous, i.e. there exists $\al:=\al(\ga)\in(0,1)$ s.t. $u\in C_{loc}^{0,\al}(\Om,\R^2).$
\label{thm:2.2.1}
\end{thm}

The second result shows, that all stationary points of the functional \eqref{eq:1.1} are of class $W_{loc}^{2,2},$ as long as the integrand remains uniformly convex. \vspace{0.5cm}
\begin{thm} Let $0< \ga<1$ and $u\in \A_{u_0}(\Om,\R^2)$ be a stationary point of the functional $I$ as defined in \eqref{eq:1.1}. Then $u\in W_{loc}^{2,2}(\Om,\R^2).$
\label{thm:2.3.DM08}
\end{thm}
From standard theory one would expect that in $2\ti2-$dimensions Meyers' Theorem would apply and one could conclude smoothness in the $0< \ga<1$ case altogether. However, this is not the case, see remark \ref{re:sec.3.1} for a detailed discussion. In both cases $0< \ga<1$  and $\ga\ge1$ one needs to rely on the higher regularity result below.\\
   
Lastly, we show that every stationary point $u$ of \eqref{eq:1.1}, which is additionally in $W_{loc}^{2,2}$ and which Jacobian is suitable Hölder continuous is already $C_{loc}^{\infty}.$ Central element is to show that $\grad^2 u$ satisfies a reverse Hölder inequality. We do so by means of measure theoretic ideas and Hardy space theory.\\
\begin{thm}[Higher order regularity]
Let $0<\ga<\infty,$ let $\Om\ss \R^2$ open, bounded, and with $\p\Om\in C^{0,1},$ let $p\in[\frac{8}{5},2)$ and the corresponding $\beta:=2(\frac{2}{p}-1)>0,$ let $u\in \A_{u_0}\cap W_{loc}^{2,2}(\Om,\R^2),$ let $x\mt d_{\grad u(x)}:=\det\grad u(x)\in C_{loc}^{0,\be}(\Om)$ and $u\in \A_{u_0}(\Om,\R^2)$ be a stationary point of the functional $I$ as defined in \eqref{eq:1.1}.\\

Then there is $\d>0$ s.t. $u\in W_{loc}^{2,2+\d}(\Om,\R^2).$ In particular, there are constants $\d=\d(\ga,\|\rho''\|_{L^\infty},p)>0$ and $C=C(\ga,\|\rho''\|_{L^\infty},p)>0$ s.t. the following reverse Hölder inequality is satisfied
\begin{equation}
\left(\fint\limits_{B_{r}}{|\grad^2 u|^{2+\d}\;dx}\right)^{\frac{1}{2+\d}}\le C\left(\fint\limits_{B_{2r}}{|\grad^2 u|^2\;dx}\right)^{\frac{p}{2}}\;\mb{for any}\; B(x,2r)\ss\ss B.
\end{equation}
Moreover, $u\in C_{loc}^{\infty}(\Om,\R^2).$\\

In particular, if $0<\ga<1$ then any stationary point $u\in\A_{u_0}$ of the functional $I$, as given in \eqref{eq:1.1}, is smooth.\\
\label{Thm:2.4.1}
\end{thm}

The assumption $x\mt d_u(x)\in C_{loc}^{0,\be}(\Om)$ can not be dropped in general. $u\in W_{loc}^{2,2}(\Om,\R^2)$ implies  $d_u(\cd)\in W_{loc}^{1,q}(\Om,\R^2)$ for any $1\le q<2,$ which does not embed into $C_{loc}^{0,\be}(\Om)$ for any $\be\in (0,1].$ However, it is possible to argue that if $u\in W_{loc}^{2,1}(\Om,\R^2)$ then $u\in C_{c}^{0}(\Om,\R^2)$ and in particular $d_u(\cd)\in C_{c}^{0}(\Om),$ see \cite[eq.(1.2)]{BF02}. See, also \cite[thm 4]{Y94}.\\

Although, there is already a rich literature concerning regularity results of poly- or the more general class of quasiconvex functionals many questions remain still open.\footnote{For the various notions of semiconvexities and their relations, see \cite{D08, R18}.} For the classical theory concerning the special case when the integrand is uniformly convex one might start with \cite{GIU03}. For regularity results concerning infinite nonlinear elasticity (that is, where the integrand $W$ depends on $d$ in such a way, that $W(d)\ra+\infty$ if $d\ra0^+$ or $d\ra+\infty$ and $W(d)=+\infty$ if $d\le 0$) see \cite{BOP91MS,Y06,JB17,BK19,FR95}. In the following we will focus only on those classes and restrict ourselves to finite elasticity.\\

The only everywhere regularity results we are aware of are:
in \cite{FH94} Fusco and Hutchinson consider a $2\ti2-$dimensional polyconvex model problem and obtain everywhere continuity for $Q-$minimzers. In \cite{CLM17} a special class of polyconvex functionals in $3\ti3-$dimension is considered and for local minimizers local boundedness is obtained. A higher order regularity result can be found in the PhD thesis \cite{C14}. If the integrand is strongly quasiconvex and $C^2,$ the boundary conditions are smooth enough and are small in some $L^p-$norm, then maximal smoothness and uniqueness can be guaranteed. Note, that our results allow more general boundary conditions ($u\in L^2$) and no smallness condition is assumed. In \cite{BOP92} incompressible condition is present higher order regularity is discussed for the special case of the Dirichtlet energy.\\

Most work for these general classes of semiconvexities are concerned with partial regularity. For some partial regularity results in certain situations with polyconvex functionals, see \cite{FH94,EM01,CY07}. The classical works regarding partial regularity for stationary points of quasiconvex functionals were obtained by Evans \cite{Ev86} and Acerbi and Fusco \cite{AF87}. Note, that the latter result actually applies to our situation and yields that all our minimizers/stationary points have to be smooth up to a nullset. Since then these results have been extended in various ways, see for example \cite{KT03, SS09}. Recently, a partial regularity result in the context of $\A-$quasiconvex functionals has been established in \cite{CG20}. Kristensen and Mingione in \cite{KM07} were additionally able to show an estimate on the dimension of the singular set, strictly less than the full dimension, for Lipschitzian $\omega-$minima for suitable quasiconvex functionals. A partial regularity result is available in the incompressible case, see \cite{EVGA99}.\\

On the other hand, using Gromov's convex integration method, Müller and Sverak in \cite{MS03} and Kristensen and Taheri \cite{KT03} constructed quasiconvex smooth integrands s.t. the (local) minimizer is everywhere Lipschitz but nowhere $C^1$. \!\!These are examples of stationary points, where the singular set has full dimension, destroying any hope on a sharper estimate on the dimension of the singular set below the space dimension for general stationary points. Szekelyhidi, Jr. \cite{Sz04} constructed an analogous counterexample with a smooth polyconvex integrand.\\

\textbf{Plan for the paper:}
In section 2 we start with the proof of theorem \ref{thm:2.2.1}. Section 3 then discusses theorem \ref{thm:2.3.DM08} and is completed with remark 3.1. In §4.1 one can find the main part of the proof of theorem \ref{Thm:2.4.1}. In §4.2 technical lemmatas essential to the proofs of the considered theorems. In §4.3 a version of the reverse Fatou's lemma is given. The section and the paper is ended with §4.4 and a description of  Hardy spaces.

\section{Hölder continuity}
As a first result we obtain local Hölder continuity. A proof is straightforward, so that it will be enough to give some details.\\

\textbf{Proof of Theorem \ref{thm:2.2.1}:}\\
Initially, the fact that test functions are allowed to be of class $W_0^{1,2}$ enables one to choose $\vp=\eta^2 (u-a),$ where $\eta\in C_c^\infty(\Om)$ denotes the standard cut-off function. Recall the definition: for $x_0\in \Om$ and $r>0$ s.t. $B(x_0,2r)\ss\ss \Om$ we want $\eta$ to be $1$ on $B(x_0,r)$ and $0$ on $\Om\setminus B(x_0,2r).$ Additionally, there exists some constant $c>0$ s.t. $|\grad \eta|\le \frac{c}{r}$ and $\supp\grad\eta \ss B(x_0,2r)\sm B(x_0,r),$ and we choose $a:=(u)_{B_{2r}\sm B_r}=\frac{1}{\La^2(B_{2r}\sm B_r)}\int\limits_{B_{2r}\sm B_r}{u\;dx}.$ Then the gradient of $\vp$ is given by $\grad \vp=2\eta \grad\eta\ot(u-a)+\eta^2\grad u.$ Plugging this into \eqref{eq:ELE1.1}, noting that $\cof A\cdot A=2d_A$ for all $A\in \R^{2\ti2}$ and collecting all $\eta^2-$ terms to one side and all terms containing $2\eta \grad\eta\ot(u-a)$ to the other side yields,
\begin{equation}
\int\limits_{\Om}{\eta^2|\grad u|^2+2\eta^2\rho'(d)d\; dx}=-\int\limits_{\Om}{2\eta(\grad u+\rho'(d)\cof \grad u) \cdot(\grad\eta\ot(u-a))\; dx}.
\label{eq:ELE1.2}
\end{equation}
The left-hand side of \eqref{eq:ELE1.2} can be estimated from below by
It is easy to see that our $\rho$ satisfies $\rho'(d)d\ge0$ for any $d\in\R.$ Hence, for the integrand of the LHS of \eqref{eq:ELE1.2} we have
\[\eta^2|\grad u|^2\le\eta^2|\grad u|^2+2\eta^2\rho'(d)d.\]
Since $\rho'$ is smooth and monotonically increasing it holds $0\le\rho'\le\gamma$ and recalling $|\grad u|=|\cof \grad u|,$ then the RHS of \eqref{eq:ELE1.2} can be controlled by
\[
C(\ga)\int\limits_{\Om}{|\eta||\grad u||u-a||\grad\eta|\; dx},
\]
where $C(\ga)=2(1+\ga)>0$ is a positive constant.\\

Together, we arrive at an inequality of the form
\[
\int\limits_{\Om}{\eta^2|\grad u|^2\; dx}\le C(\ga) \int\limits_{\Om}{|\eta||\grad u||u-a||\grad\eta|\; dx}.
\]
The proof is completed by following the standard Caccioppoli method. \vspace{1cm}

\section{From $W^{1,2}$ to $W_{loc}^{2,2}$ for any $0<\ga<1$ via difference quotients}
\label{sec:2.3}

In the last paragraph we have seen that all stationary points must be locally Hölder continuous. It is natural to ask whether the regularity can be improved even further. In this section we show that if $0<\ga<1,$ i.e. the integrand is uniformly convex, then stationary points of the functional $I,$ as defined in \eqref{eq:1.1}, are in $W_{loc}^{2,2}$. This agrees from what we would expect by standard theory. But one of the requirements to apply that theory is that $\grad_\xi^2W(\grad u(\cd))$ needs to be locally bounded (in $L^\infty$) from above. But $\xi\mapsto\grad_\xi^2W(\xi)$ rather grows quadratically in $\xi.$\\

Why does the standard method fail precisely? Usually, one would like to test the ELE \eqref{eq:ELE1.1} with $\vp^s(x):=-\D^{-h,s}(\eta^2(x)\D^{h,s}u(x)),$ where $\D^{h,s}f(x)=h^{-1}(f(x+he_s)-f(x))$ and $\D^{-h,s}f(x)=h^{-1}(f(x)-f(x-he_s))$ for a.e. $x\in\Om_h:=\{x\in\Om:\dist(x,\p\Om)>h\},$ $h\in \R_0^+,$ $s=1,2,$ and where $\{e_1,e_2\}$ denotes the standard basis. We can then simply derive the following estimate
\begin{eqnarray}
(1-\ga)\int\limits_{\Om_h}{\eta^2|\grad \Delta^{h,s}u|^2\; dx}\le C\int\limits_{\Om_h}{|\eta||\Delta^{h,s}\grad_\xi W(\grad u)||\grad\eta||\Delta^{h,s}u| \; dx},
\label{eq:2.3.1001}
\end{eqnarray}
which is ok if $\ga\in(0,1).$ Then in order for the RHS to be bounded, we need $\Delta^{h,s}\grad_\xi W(\grad u)\in L^2(\Om_h),$ which can be guaranteed if $D_\xi^2 W(\grad u(\cdot))\in L^\infty(\Om_h).$ Unfortunately, since $|D_\xi^2 W(\grad u(\cd))|\sim|\grad u(\cd)|^2,$ $D_\xi^2 W(\grad u(\cdot))$ is in $L^1(\Om_h),$ but not necessarily in any better space. Hence, the RHS of \eqref{eq:2.3.1001} might blow up.\\

Situations, where the second derivative of a uniformly convex integrand might be unbounded, have been studied in \cite{DM09}. De Maria's genius idea is to play the problematic term against a slight variation of itself. In order to explain the idea even further consider \begin{equation}\int\limits_{\Om_h}{(\eta^2(x+he_s)-\eta^2(x))f(\D^{h,s}D_\xi W(\grad u(x)))\;dx}, \label{eq:2.3.1000a}\end{equation}
$f$ denoting a generic function. 
Now for the most part of the set $\Om_h$ the integrand in \eqref{eq:2.3.1000a} disappears. Even though the integrand $f$ might be large, there is only a small set on which $f$ can contribute. As a consequence, this introduces another `smallness' to the integral, which turns out to be enough to control its behaviour.\\
 
As a test function we will therefore use 
\begin{equation}\vp(x)=-\D^{-h,s}((\eta^2(x+he_s)-(1-\al)\eta^2(x))\D^{h,s} u(x)),\label{eq:2.3.1000b}\end{equation}
a.e. in $\Om_h$ and for some $0<\al<1.$ Notice, that we use a slight modification of the test function given in \cite{DM09}. The reason is, that we do not explicitly use the strong ellipticity inequality, which can be avoided by using the concrete form of our integrand $W$ and by introducing an additional small parameter $\al.$ The advantage being, that our method could be used in future works when $\ga\ge 1.$ But, at present, this remains open. We want to emphasise another time, that the condition $0<\ga<1$ will be crucial in the following discussion.\\

The goal is to obtain a reverse Poincar\capitalacute e inequality of the form
\[\|\grad^2 u\|_{L^2(B_r)}^2\le \frac{C(\ga)}{r^2}\|\grad u\|_{L^2(B_{2r})}^2 \;\mb{for any}\; B(2r)\ss\ss \Om.\]

The strategy of the proof is then fairly straightforward: 
\begin{enumerate}
\item Testing the ELE with \eqref{eq:2.3.1000b} and organising the upcoming terms. 
\item Estimating, mostly via standard estimates, however, some of the estimates are more subtle, where explicit use of the closeness of $\eta(x)$ and $\eta(x+he_s)$ must be made, thereby following closely \cite{DM09}.
\item The proof is then concluded by Nirenberg's lemma.\\
\end{enumerate}

\textbf{Proof of theorem \ref{thm:2.3.DM08}:}\\
\textbf{Step 1: testing the equation.}
Recall the ELE
\begin{equation}
\int\limits_{\Om}{(\grad u+\rho'(d)\cof \grad u) \cdot\grad\vp\; dx}=0\; {\mbox{for all}}\; \vp\in W_0^{1,2}(\Om,\R^2).
\end{equation}

Let $x_0\in \Om$ and $r>0$ s.t. $B(x_0,3r)\ss\ss \Om.$ We will use the notation $\xb:=x+he_s$ and $x^-:=x-he_s$ for every $x\in\Om_h:=\{x\in\Om:\dist(x,\p\Om)>h\}.$ Further, choose $h_0:=\frac{r}{10}$ and let $0<h<h_0,$ then by construction $\supp{\eta(\cdot)},\supp{\eta(\cd\pm he_s)}\ss\ss B(x_0,\frac{7r}{4}).$ As a test function, motivated by \cite{DM09}, we choose 
\begin{equation}\vp=-\D^{-h,s}(\ta^h(x)\D^{h,s} u),\;\mb{with}\; \ta^h(x):=\eta^2(\xb)-(1-\al)\eta^2(x) \label{eq:2.3.98}\end{equation}
where $0<\al<1$ and $\eta\in C_c^\infty(\Om)$ is a standard mollifier satisfying the properties $\eta\equiv1$ in $B_r=B(x_0,r),$ $\eta\equiv0$ in $\Om\sm B(x_0,\frac{3r}{2}),$ $0\le\eta\le1$ and there exists $c>0$ s.t. $|\grad\eta|\le\frac{c}{r}$ and $|\grad^2\eta|\le\frac{c}{r^2}.$
Differentiating \eqref{eq:2.3.98} gives
\[\grad\vp=-\D^{-h,s}(\ta^h\D^{h,s} \grad u)-\D^{-h,s}(\grad \ta^h\ot\D^{h,s} u),\]

Then the ELE becomes{\small{
\begin{equation}{
\small{-\int\limits_{B_{2r}}{(\grad u+\rho'(d)\cof \grad u)\cd\D^{-h,s}(\ta^h\D^{h,s} \grad u)\;dx}
=\int\limits_{B_{2r}}{(\grad u+\rho'(d)\cof \grad u)\cd\D^{-h,s}(\grad \ta^h\ot\D^{h,s} u)\;dx}.}}
\label{eq:2.3.99}
\end{equation}}}
First recall the integration by parts formula for difference quotients for two functions $v,w\in L^2(\Om)$ of which at least one has support in $\Om_h,$ then 
\[\int\limits_{\Om}{v\D^{h,s}w\;dx}=-\int\limits_{\Om}{(\D^{-h,s}v)w\;dx}.\]
Applying this to the LHS of \eqref{eq:2.3.99} gives
\begin{flalign*}
LHS=&\int\limits_{B_{2r}}{\ta^{h}\D^{h,s}(\grad u+\rho'(d)\cof \grad u)\cd\D^{h,s}\grad u\;dx}&\\
=&\int\limits_{B_{2r}}{\eta^2(\xb)\D^{h,s}\grad u(x)\cd\D^{h,s}\grad u\;dx}-(1-\al)\int\limits_{B_{2r}}{\eta^2(x)\D^{h,s}\grad u(x)\cd\D^{h,s}\grad u\;dx}&\\
&+\int\limits_{B_{2r}}{\ta^{h}\D^{h,s}(\rho'(d)\cof \grad u)(x)\cd\D^{h,s}\grad u\;dx},&
\end{flalign*}
where we used the definition of $\ta^h.$ 

Recalling the product rule for difference quotients, for $f:\Om_h\rightarrow \R$ and $g:\Om\rightarrow \R$ it holds that
\begin{equation}
\D^{h,s}(fg)(x)=\D^{h,s}(f)(x)g(x)+f(\bar{x})\D^{h,s}g(x).
\label{eq:2.3.100PRFD}
\end{equation}
We point out, that this formula remains true even if $f$ and $g$ are vector- or matrixvalued and for various types of products (for example, scalar- or tensorproducts). Applying the product rule to the RHS of \eqref{eq:2.3.99} yields
\begin{flalign*}
RHS=&\int\limits_{B_{2r}}{(\grad u+\rho'(d)\cof \grad u)\cd\grad\ta^{h}(x)\ot\D^{-h,s}\D^{h,s} u(x)\;dx}&\\
&+\int\limits_{B_{2r}}{(\grad u+\rho'(d)\cof \grad u)\cd\D^{-h,s}(\grad\ta^{h})(x)\ot\D^{h,s} u(x^-)\;dx}.&
\end{flalign*}
By plugging this back into the ELE \eqref{eq:2.3.99} and rearranging terms we get
\begin{flalign*}
\int\limits_{B_{2r}}{\eta^2(\xb)|\D^{h,s}\grad u|^2(x)\;dx}=&(1-\al)\int\limits_{B_{2r}}{\eta^2(x)|\D^{h,s}\grad u|^2(x)\;dx}&\\
&-\int\limits_{B_{2r}}{\ta^{h}\D^{h,s}(\rho'(d)\cof \grad u)(x)\cd\D^{h,s}\grad u\;dx}&\\
&+\int\limits_{B_{2r}}{(\grad u+\rho'(d)\cof \grad u)\cd\grad\ta^{h}(x)\ot\D^{-h,s}\D^{h,s} u(x)\;dx}&\\
&+\int\limits_{B_{2r}}{(\grad u+\rho'(d)\cof \grad u)\cd\D^{-h,s}(\grad\ta^{h})(x)\ot\D^{h,s} u(x^-)\;dx}&\\
=:&(i)+(ii)+(iii)+(iv).&
\end{flalign*}
\textbf{Step 2: about $(i).$}
We start with the following estimate: for $\be\in\R,$ we have
\begin{flalign*}\be^2\eta^2(x)=|\be(\eta(x)-\eta(\xb))+\be\eta(\xb)|^2
=\be^2(\eta(x)-\eta(\xb))^2+\be^2\eta^2(\xb)+2\be^2(\eta(x)-\eta(\xb))\eta(\xb).
\end{flalign*}
Then by Cauchy's inequality we get, for some $c'>0,$
\begin{flalign*}2\be^2(\eta(x)-\eta(\xb))\eta(\xb)\le\frac{1}{c'}(\eta(x)-\eta(\xb))^2+c'\be^4\eta^2(\xb).
\end{flalign*}
Choosing $\be^2=(1-\al)$ yields,
\begin{flalign*}(1-\al)\eta^2(x)\le\left((1-\al)+\frac{1}{c'}\right)(\eta(x)-\eta(\xb))^2+(1-\al)(1+c'(1-\al))\eta^2(\xb).
\end{flalign*}

Then by the latter estimate we can estimate $(i)$ by
\begin{flalign}
(1-\al)\int\limits_{B_{2r}}{\eta^2(x)|\D^{h,s}\grad u|^2\;dx}
\le& \left((1-\al)+\frac{1}{c'}\right)\int\limits_{B_{2r}}{|\eta(\xb)-\eta(x)|^2|\D^{h,s}\grad u|^2\;dx}&\nonumber\\
&+(1-\al)(1+c'(1-\al))\int\limits_{B_{2r}}{\eta(\xb)^2|\D^{h,s}\grad u|^2\;dx}&\nonumber\\
=:&(i.a)+(i.b).&\label{eq:2.3.1000a+b}
\end{flalign}
Using the mean value theorem and $|\grad \eta|\le\frac{c}{r}$ we can control $(i.a)$ via
\begin{flalign}
(i.a)\le&h^2\left((1-\al)+\frac{1}{c'}\right)\int\limits_{B_{\frac{7r}{4}}}{|\grad\eta|^2|\D^{h,s}\grad u|^2\;dx}&\nonumber\\
\le&\frac{h^2c^2\left((1-\al)+\frac{1}{c'}\right)}{r^2}\int\limits_{B_{\frac{7r}{4}}}{|\D^{h,s}\grad u|^2\;dx}&\nonumber\\
\le&\frac{c^2\left((1-\al)+\frac{1}{c'}\right)}{r^2}\int\limits_{B_{\frac{7r}{4}}}{|\grad u(\xb)|^2+|\grad u(x)|^2\;dx}&\nonumber\\
\le&\frac{c^2\left((1-\al)+\frac{1}{c'}\right)}{r^2}\int\limits_{B_{2r}}{|\grad u(x)|^2\;dx},&
\label{eq:2.3.1000.a.2}
\end{flalign}
where we used
\[\int\limits_{B_{\frac{7r}{4}}}{|\grad u(\xb)|^2\;dx}\le\int\limits_{B_{2r}}{|\grad u(x)|^2\;dx}.\]
From \eqref{eq:2.3.1000a+b} together with \eqref{eq:2.3.1000.a.2} we get
\begin{equation}
(i)\le (1-\al)(1+c')\int\limits_{B_{2r}}{\eta^2(\xb)|\D^{h,s}\grad u|^2\;dx}+\frac{C}{r^2}\int\limits_{B_{2r}}{|\grad u(x)|^2\;dx},
\label{eq:2.3.100.i}
\end{equation}
where we simplified the coefficient in front of the first integral by using the estimate $(1+c'(1-\al))\le(1+c').$ It will be crucial that for every $0<\al<1$ it is possible to choose $c'>0$ so small that $(1-\al)(1+c')<1.$ Moreover, we introduced a generic constant $C>0,$ which is also allowed to depend on $\ga.$ We can be generous, with the constant in front of the rightmost integral, as long as it remains finite.\\

\textbf{Step 3: about $(ii).$}
Rewriting
$\ta^h(x)=\al\eta^2(\xb)+(1-\al)(\eta^2(\xb)-\eta^2(x))$
yields
\begin{flalign*}
(ii)=&-\int\limits_{B_{2r}}{\ta^h(x)\D^{h,s}(\rho'(d)\cof \grad u)(x)\cd\D^{h,s}\grad u\;dx}&\\
=&-\al\int\limits_{B_{2r}}{\eta^2(\xb)\D^{h,s}(\rho'(d)\cof \grad u)(x)\cd\D^{h,s}\grad u\;dx}&\\
&-(1-\al)\int\limits_{B_{2r}}{(\eta^2(\xb)-\eta^2(x))\D^{h,s}(\rho'(d)\cof \grad u)(x)\cd\D^{h,s}\grad u\;dx}.&
\end{flalign*}

For the first term we use the lower bound of lemma \ref{HR:LUBs}, which is given by
\begin{equation}
\D^{h,s}(\rho'(d)\cof \grad u)(x)\cd\D^{h,s}\grad u(x)\ge-\ga|\D^{h,s}\grad u(x)|^2\end{equation}
for the second one we use the upper bound, shown in the same lemma, given by
\begin{equation}
|\Delta^{h,s}(\rho'(d)\cof\grad u)(x)|\le\ga|\Delta^{h,s}\grad u(x)|+\frac{2\ga}{h}|\grad u(x)|
\label{eq:2.16}
\end{equation}
which hold a.e. in $\Om_{h}.$ Hence,
\begin{flalign*}
(ii)\le&\al\ga\int\limits_{B_{2r}}{\eta^2(\xb)|\D^{h,s}\grad u|^2\;dx}&\\
&+(1-\al)\ga\int\limits_{B_{2r}}{|\eta^2(\xb)-\eta^2(x)||\D^{h,s}\grad u|^2\;dx}&\\&+\frac{2(1-\al)\ga}{h}\int\limits_{B_{2r}}{|\eta^2(\xb)-\eta^2(x)||\grad u(x)||\D^{h,s}\grad u|\;dx}&\\
=:&(ii.a)+(ii.b)+(ii.c)&
\end{flalign*}
(ii.a) is already of the desired form. For (ii.b) we apply, again, the mean value theorem and $|\grad\eta|\le \frac{c}{r}$ to obtain

\begin{flalign}
(ii.b)\le&(1-\al)\ga\int\limits_{B_{2r}}{|\eta(x)-\eta(\xb)||\eta(x)+\eta(\xb)||\D^{h,s}\grad u|^2\;dx}&\nonumber\\
\le&(1-\al)\ga h\int\limits_{B_{\frac{7r}{4}}}{|\grad\eta||\eta(x)+\eta(\xb)||\D^{h,s}\grad u|^2\;dx}&\nonumber\\
\le&Ch^2\ga^2\int\limits_{B_{\frac{7r}{4}}}{|\grad\eta||\D^{h,s}\grad u(x)|^2\;dx}+\ve(1-\al)^2\int\limits_{B_{2r}}{|\eta(x)+\eta(\xb)|^2|\D^{h,s}\grad u(x)|^2\;dx}&\nonumber\\
\le&C(\ga) h^2\int\limits_{B_{\frac{7r}{4}}}{|\grad\eta|^2|\D^{h,s}\grad u(x)|^2\;dx}+2\ve(1-\al)^2\int\limits_{B_{2r}}{(\eta^2(x)+\eta^2(\xb))|\D^{h,s}\grad u(x)|^2\;dx}.&\label{eq:2.16.a}
\end{flalign}
Here the $\ve$ results from another application of Cauchy's inequality with the weight $\ve>0,$ where we can choose $\ve$ arbitrarily small. The counterpart $\ve^{-1}$ has been absorbed into $C.$\\
Then the first term of the latter line can be controlled, as before in \eqref{eq:2.3.1000.a.2}, by
\begin{flalign*}
C(\ga) h^2\int\limits_{B_{\frac{7r}{4}}}{|\grad\eta|^2|\D^{h,s}\grad u(x)|^2\;dx}\le \frac{C(\ga)}{r^2} \int\limits_{B_{\frac{7r}{4}}}{|\grad u(\xb)|^2+|\grad u(x)|^2\;dx}\le\frac{C(\ga)}{r^2} \int\limits_{B_{2r}}{|\grad u(x)|^2\;dx}.
\end{flalign*}
The latter term of \eqref{eq:2.16.a} can be estimated in a similar fashion as has been done for $(i):$
\begin{flalign*}
&2\ve(1-\al)^2\int\limits_{B_{2r}}{(\eta^2(x)+\eta^2(\xb))|\D^{h,s}\grad u(x)|^2\;dx}&\\
\le&(2\ve(1-\al)^2+2\ve(1-\al)^2(1+2c'\ve(1-\al)^2))\int\limits_{B_{2r}}{\eta^2(\xb)|\D^{h,s}\grad u|^2\;dx}+\frac{c}{r^2}\int\limits_{B_{2r}}{|\grad u(x)|^2\;dx}.&
\end{flalign*}
Together we get
\begin{equation}
(ii.b)\le (2\ve(1-\al)^2(2+2c''\ve(1-\al)^2))\int\limits_{B_{2r}}{\eta^2(\xb)|\D^{h,s}\grad u(x)|^2\;dx}+\frac{C(\ga)}{r^2} \int\limits_{B_{2r}}{|\grad u(x)|^2\;dx}.
\label{eq:2.3.100}
\end{equation}

A similar calculation can be done for (ii.c):
\begin{flalign}
(ii.c)=&+\frac{2(1-\al)\ga}{h}\int\limits_{B_{2r}}{|\eta^2(\xb)-\eta^2(x)||\grad u(x)||\D^{h,s}\grad u|\;dx}&\nonumber\\
\le&2(1-\al)\ga \int\limits_{B_{\frac{7r}{4}}}{|\grad\eta||\eta(x)+\eta(\xb)||\grad u(x)||\D^{h,s}\grad u|\;dx}&\nonumber\\
\le&C\ga^2\int\limits_{B_{\frac{7r}{4}}}{|\grad\eta|^2|\grad u(x)|^2\;dx}+\ve(1-\al)^2\int\limits_{B_{2r}}{|\eta(x)+\eta(\xb)|^2|\D^{h,s}\grad u(x)|^2\;dx}&\nonumber\\
\le&(2\ve(1-\al)^2(2+2c''\ve(1-\al)^2))\int\limits_{B_{2r}}{\eta^2(\xb)|\D^{h,s}\grad u(x)|^2\;dx}+\frac{C(\ga)}{r^2} \int\limits_{B_{2r}}{|\grad u(x)|^2\;dx}.&
\label{eq:2.3.101}
\end{flalign}

Collecting, (ii.a), \eqref{eq:2.3.100} and \eqref{eq:2.3.101} yields, 
\begin{equation}
(ii)\le (\al\ga+4\ve(1-\al)^2(2+2c''\ve(1-\al)^2))\int\limits_{B_{2r}}{\eta^2(\xb)|\D^{h,s}\grad u(x)|^2\;dx}+\frac{C(\ga)}{r^2} \int\limits_{B_{2r}}{|\grad u(x)|^2\;dx}.
\label{eq:2.3.100.ii}
\end{equation}

\textbf{Step 4: about $(iii).$} The third term can be estimated by
\begin{flalign*}
(iii)=&\int\limits_{B_{2r}}{(\grad u+\rho'(d)\cof \grad u)\cd\grad\ta^{h}(x)\ot\D^{-h,s}\D^{h,s} u(x)\;dx}&\\
\le&(1+\ga)\int\limits_{B_{2r}}{|\grad\ta^{h}(x)||\grad u||\D^{-h,s}\D^{h,s} u|\;dx}.&
\end{flalign*}
Next $\grad \ta^h$ is given by
\[
\grad\ta^{h}(x)=2\eta(\xb)\grad\eta(\xb)-2(1-\al)\eta(x)\grad\eta(x), \;{\mb{for all}}\; x\in\Om_h.
\]
By the latter expression and the triangle inequality we get
\begin{flalign*}
(iii)\le&2(1+\ga)\int\limits_{B_{2r}}{\eta(\xb)|\grad\eta(\xb)||\grad u||\D^{-h,s}\D^{h,s} u|\;dx}&\\
&+2(1+\ga)(1-\al)\int\limits_{B_{2r}}{\eta(x)|\grad\eta(x)||\grad u||\D^{-h,s}\D^{h,s} u|\;dx}&\\
=:&(iii.a)+(iii.b)&
\end{flalign*}
Both terms are very similar, let's start with the second one. By Young's inequality and $|\grad\eta|\le\frac{c}{r}$ we get
\begin{flalign*}
(iii.b)\le&\ve\int\limits_{B_{2r}}{\eta^2(x)|\D^{-h,s}\D^{h,s} u|^2\;dx}+\frac{(1+\ga)^2(1-\al)^2}{\ve}\int\limits_{B_{2r}}{|\grad\eta(x)|^2|\grad u(x)|^2\;dx}&\\
\le&\ve\int\limits_{B_{2r}}{\eta^2(x)|\D^{-h,s}\D^{h,s} u|^2\;dx}+\frac{c(1+\ga)^2}{\ve r^2}\int\limits_{B_{2r}}{|\grad u(x)|^2\;dx}&
\end{flalign*}
The first term on the RHS seems difficult. But to our advantage De Maria shows in \cite{DM09} how this term can be controlled. First we rewrite the integrand by
\begin{flalign}
&\ve\int\limits_{B_{2r}}{\eta^2(x)|\D^{-h,s}\D^{h,s} u(x)|^2\;dx}&\nonumber\\
=&\ve\int\limits_{B_{2r}}{\frac{1}{h^2}|\eta(x)\D^{h,s} u(x)-\eta(x)\D^{h,s} u(x^-)|^2\;dx}&\nonumber\\
=&\ve\int\limits_{B_{2r}}{\frac{1}{h^2}|\eta(x^-)\D^{h,s}u(x^-)-\eta(x)\D^{h,s}u(x)+(\eta(x)-\eta(x^-))\D^{h,s} u(x^-)|^2\;dx}&\nonumber\\
\le&2\ve\int\limits_{B_{2r}}{|\D^{-h,s}(\eta\D^{h,s}u)(x)|^2\;dx}+2\ve\int\limits_{B_{2r}}{|\D^{-h,s}\eta(x)|^2|\D^{h,s}u(x^-)|^2\;dx},&\label{HR:(iii).1}
\end{flalign}
where we used $(a+b)^2\le2(a^2+b^2)$ and the second term of \eqref{HR:(iii).1} can be covered as follows
\begin{flalign*}
2\ve\int\limits_{B_{2r}}{|\D^{-h,s}\eta(x)|^2|\D^{h,s}u(x^-)|^2\;dx}\le&\frac{2C\ve}{r^2}\int\limits_{B_{\frac{7r}{4}}}{|\D^{h,s}u(x^-)|^2\;dx}&\\
\le&\frac{2C\ve}{r^2}\int\limits_{B_{\frac{7r}{4}+h_0}}{|\D^{h,s}u(x)|^2\;dx}&\\
\le&\frac{2C\ve}{r^2}\int\limits_{B_{2r}}{|\grad u(x)|^2\;dx},&
\end{flalign*}
where we used Nirenberg's lemma\footnote{see \cite[lem 7.23-7.24]{GT}.} in the latter estimate. The first term of \eqref{HR:(iii).1} can be treated, again, using Nirenberg's lemma and the product rule, as follows
\begin{flalign*}
2\ve\int\limits_{B_{2r}}{|\D^{-h,s}(\eta\D^{h,s}u)(x)|^2\;dx}
\le&2\ve\int\limits_{B_{2r}}{|\grad(\eta\D^{h,s}u)(x)|^2\;dx}&\\
\le&4\ve\int\limits_{B_{2r}}{|\grad\eta(x)\D^{h,s}u(x)|^2\;dx}+4\ve\int\limits_{B_{2r}}{\eta^2(x)|\grad\D^{h,s}u(x)|^2\;dx}&\\
\le&\frac{c\ve}{r^2}\int\limits_{B_{\frac{7r}{4}}}{|\D^{h,s}u(x)|^2\;dx}+4\ve\int\limits_{B_{2r}}{\eta^2(x)|\grad\D^{h,s}u(x)|^2\;dx}&\\
\le&\frac{c\ve}{r^2}\int\limits_{B_{2r}}{|\grad u(x)|^2\;dx}+c\ve\int\limits_{B_{2r}}{\eta^2(\xb)|\grad\D^{h,s}u(x)|^2\;dx}.&
\end{flalign*}
In the last step we used Nirenberg to bound the first term, and step (i) in the rightmost term. Hence,
\begin{flalign*}
(iii.b)\le c\ve\int\limits_{B_{2r}}{\eta^2(\xb)|\D^{h,s}\grad u|^2\;dx}+\frac{C(\ga)}{r^2}\int\limits_{B_{2r}}{|\grad u(x)|^2\;dx}.
\end{flalign*}
The (iii.a)-term can be controlled analogously:
\begin{flalign*}
(iii.a)\le&\ve\int\limits_{B_{2r}}{\eta^2(\xb)|\D^{-h,s}\D^{h,s} u|^2\;dx}+\frac{(1+\ga)^2}{\ve}\int\limits_{B_{2r}}{|\grad\eta(\xb)|^2|\grad u(x)|^2\;dx}&\\
\le&\ve\int\limits_{B_{2r}}{\eta^2(\xb)|\D^{-h,s}\D^{h,s} u|^2\;dx}+\frac{c(1+\ga)^2}{\ve r^2}\int\limits_{B_{2r}}{|\grad u(x)|^2\;dx}&
\end{flalign*}
The first term can be covered by
\begin{flalign}
&\ve\int\limits_{B_{2r}}{\eta^2(\xb)|\D^{-h,s}\D^{h,s} u(x)|^2\;dx}&\nonumber\\
=&\ve\int\limits_{B_{2r}}{\frac{1}{h^2}|\eta(\xb)\D^{h,s} u(x)-\eta(\xb)\D^{h,s} u(x^-)|^2\;dx}&\nonumber\\
=&\ve\int\limits_{B_{2r}}{\frac{1}{h^2}|\eta(x)\D^{h,s}u(x^-)-\eta(\xb)\D^{h,s}u(x)+(\eta(\xb)-\eta(x))\D^{h,s} u(x^-)|^2\;dx}&\nonumber\\
\le&2\ve\int\limits_{B_{2r}}{|\D^{h,s}(\eta(x)\D^{h,s}u(x-he_s))|^2\;dx}+2\ve\int\limits_{B_{2r}}{|\D^{h,s}\eta(x)\D^{h,s}u(x^-)|^2\;dx}&\label{eq:2.3.102}
\end{flalign}
First term of \eqref{eq:2.3.102}:
\begin{flalign*}
&2\ve\int\limits_{B_{2r}}{|\D^{h,s}(\eta(x)\D^{h,s}u(x-he_s))|^2\;dx}&\\
\le&2\ve\int\limits_{B_{2r}}{|\grad(\eta(x)\D^{h,s}u(x^-))|^2\;dx}&\\
\le&4\ve\int\limits_{B_{2r}}{|\grad\eta(x)\D^{h,s}u(x^-)|^2\;dx}+4\ve\int\limits_{B_{2r}}{\eta^2(x)|\grad\D^{h,s}u(x^-)|^2\;dx}&\\
\le&\frac{c\ve}{r^2}\int\limits_{B_{2r}}{|\grad u(x)|^2\;dx}+c\ve\int\limits_{B_{2r}}{\eta^2(\xb)|\grad\D^{h,s}u(x)|^2\;dx}&
\end{flalign*}
Rightmost term of \eqref{eq:2.3.102}:
\begin{flalign*}
2\ve\int\limits_{B_{2r}}{|\D^{h,s}\eta(x)|^2|\D^{h,s}u(x^-)|^2\;dx}\le&\frac{c\ve}{r^2}\int\limits_{B_{\frac{7r}{4}}}{|\D^{h,s}u(x^-)|^2\;dx}&\\
\le&\frac{c}{r^2}\int\limits_{B_{\frac{7r}{4}+h_0}}{|\D^{h,s}u(x)|^2\;dx}&\\
\le&\frac{c}{r^2}\int\limits_{B_{2r}}{|\grad u(x)|^2\;dx}&
\end{flalign*}

Hence,
\[
(iii.a)\le c\ve\int\limits_{B_{2r}}{\eta^2(\xb)|\D^{h,s}\grad u|^2\;dx}+\frac{C(\ga)}{r^2}\int\limits_{B_{2r}}{|\grad u(x)|^2\;dx}.
\]

Combining the bounds for (iii.a) and (iii.b) yields
\begin{equation}
(iii)\le c\ve\int\limits_{B_{2r}}{\eta^2(\xb)|\D^{h,s}\grad u(x)|^2\;dx}+\frac{C(\ga)}{r^2} \int\limits_{B_{2r}}{|\grad u(x)|^2\;dx}.
\label{eq:2.3.100.iii}
\end{equation}

\textbf{Step 5: about $(iv).$} Starting, similarly to $(iii)$ by
\begin{flalign*}
(iv)=&\int\limits_{B_{2r}}{(\grad u+\rho'(d)\cof \grad u)\cd\D^{-h,s}(\grad\ta^{h})(x)\ot\D^{h,s} u(x^-)\;dx}&\\
\le&(1+\ga)\int\limits_{B_{2r}}{|\grad u||\D^{-h,s}(\grad\ta^{h})(x)||\D^{h,s} u(x^-)|\;dx}.&
\end{flalign*}
We can write
\begin{flalign*}
\D^{-h,s}(\grad\ta^{h})(x)=&2\D^{-h,s}[\eta(\xb)\grad\eta(\xb)-(1-\al)\eta(x)\grad\eta(x)]&\\
=&2[\D^{-h,s}(\eta(\xb))\grad\eta(\xb)+\eta(\xb)\D^{-h,s}(\grad\eta(\xb))&\\
&-(1-\al)\D^{-h,s}(\eta(x))\grad\eta(x)-(1-\al)\eta(x)\D^{-h,s}(\grad\eta(x))].&
\end{flalign*}
Using this latter expansion we can estimate $(iv)$ by
\begin{flalign*}
(iv)=&\int\limits_{B_{2r}}{(\grad u+\rho'(d)\cof \grad u)\cd\D^{-h,s}(\grad\ta^{h})(x)\ot\D^{h,s} u(x^-)\;dx}&\\
\le&2(1+\ga)\int\limits_{B_{2r}}{|\grad u||\D^{-h,s}(\eta(\xb))||\grad\eta(x)||\D^{h,s} u(x^-)|\;dx}&\\
&+2(1+\ga)\int\limits_{B_{2r}}{|\grad u||\eta(\xb)||\D^{-h,s}(\grad\eta(\xb))||\D^{h,s} u(x^-)|\;dx}&\\
&+2(1-\al)(1+\ga)\int\limits_{B_{2r}}{|\grad u||\D^{-h,s}(\eta(x))||\grad\eta(x)||\D^{h,s} u(x^-)|\;dx}&\\
&+2(1-\al)(1+\ga)\int\limits_{B_{2r}}{|\grad u||\eta(x^-)||\D^{-h,s}(\grad\eta(x))||\D^{h,s} u(x^-)|\;dx}&\\
=:&(iv.a)+(iv.b)+(iv.c)+(iv.d).&
\end{flalign*}
We can treat the first part in the following way 
\begin{flalign*}
(iv.a)=&2(1+\ga)\int\limits_{B_{2r}}{|\grad u||\D^{-h,s}(\eta(\xb))||\grad\eta(x)||\D^{h,s} u(x^-)|\;dx}&\\
\le&C(\ga)\int\limits_{B_{2r}}{|\grad\eta(x)|^2|\grad u|^2\;dx}+\ve\int\limits_{B_{2r}}{|\D^{-h,s}(\eta(\xb))|^2|\D^{h,s} u(x^-)|^2\;dx}&\\
\le&\frac{C(\ga)}{r^2}\int\limits_{B_{2r}}{|\grad u|^2\;dx}+\ve\int\limits_{B_{\frac{7r}{4}}}{|\grad\eta|^2|\D^{h,s} u(x^-)|^2\;dx}&\\
\le&\frac{C(\ga)}{r^2}\int\limits_{B_{2r}}{|\grad u|^2\;dx}+\frac{\ve}{r^2}\int\limits_{B_{\frac{7r}{4}}}{|\D^{h,s} u(x^-)|^2\;dx}&\\
\le&\frac{C(\ga)}{r^2}\int\limits_{B_{2r}}{|\grad u|^2\;dx}+\frac{c\ve}{r^2}\int\limits_{B_{2r}}{|\grad u(x)|^2\;dx}&\\
\le&\frac{C(\ga)}{r^2}\int\limits_{B_{2r}}{|\grad u|^2\;dx},&
\end{flalign*}
where we used Young's inequality, the mean value theorem $|\grad\eta|<\frac{c}{r}$ and Nirenberg's lemma.
Similarly we get \vspace{-1mm}
\begin{align*}
(iv.c)\le&\frac{C(\ga)}{r^2}\int\limits_{B_{2r}}{|\grad u|^2\;dx}.&
\end{align*}
For (iv.b) we use Young's inequality and the mean value theorem but this time we apply it to $\grad \eta,$ combined with $|\grad^2\eta|\le \frac{c}{r^2}$ we get
\begin{flalign*}
(iv.b)=&2(1+\ga)\int\limits_{B_{2r}}{|\grad u||\eta(\xb)||\D^{-h,s}(\grad\eta(\xb))||\D^{h,s} u(x^-)|\;dx}&\\
\le&C(\ga)\int\limits_{B_{2r}}{|\D^{-h,s}(\grad\eta(\xb))||\grad u|^2\;dx}+\ve\int\limits_{B_{2r}}{|\eta(\xb)||\D^{-h,s}(\grad\eta(\xb))||\D^{h,s} u(x^-)|^2\;dx}&\\
\le&\frac{C(\ga)}{r^2}\int\limits_{B_{2r}}{|\grad u|^2\;dx}+\frac{c\ve}{r^2}\int\limits_{B_{2r}}{|\eta(\xb)|^2|\D^{h,s} u(x^-)|^2\;dx}&\\
\le&\frac{C(\ga)}{r^2}\int\limits_{B_{2r}}{|\grad u|^2\;dx}+\frac{c\ve}{r^2}\int\limits_{B_{\frac{7r}{4}}}{|\D^{h,s} u(x^-)|^2\;dx}&\\
\le&\frac{C(\ga)}{r^2}\int\limits_{B_{2r}}{|\grad u|^2\;dx}&
\end{flalign*}
Similarly, we obtain
\begin{align*}
(iv.d)\le&\frac{C(\ga)}{r^2}\int\limits_{B_{2r}}{|\grad u|^2\;dx}.&
\end{align*}

Together,\vspace{-3mm}
\begin{equation}
(iv)\le \frac{C(\ga)}{r^2} \int\limits_{B_{2r}}{|\grad u(x)|^2\;dx}.
\label{eq:2.3.100.iv}
\end{equation}
\textbf{Step 6: conclusion.} Collecting, \eqref{eq:2.3.100.i}, \eqref{eq:2.3.100.ii}, \eqref{eq:2.3.100.iii}, and \eqref{eq:2.3.100.iv} yields the following inequality
\begin{flalign}
\int\limits_{B_{2r}}{\eta^2(\xb)|\D^{h,s}\grad u|^2(x)\;dx}\le\wt{C}\int\limits_{B_{2r}}{\eta^2(\xb)|\D^{h,s}\grad u(x)|^2\;dx}
+\frac{C(\ga)}{r^2} \int\limits_{B_{2r}}{|\grad u(x)|^2\;dx},\label{eq:2.3.100.v}
\end{flalign}
where $\wt{C}:=\wt{C}(\al,c',c,\ga,\ve):=((1-\al)(1+c')+\al\ga+c\ve)>0.$
Realising that $\wt{C}<1$ if $0\le\ga<1$ and $\al,\ve$ small enough\footnote{This can be seen if for any $0<\al,\ga<1$ we choose $c'>0$ s.t. $1+c'=(1-\al(1-\ga))^{-1}.$ Then 
\[1-((1-\al)(1+c')+\al\ga)=\frac{\al^2\ga(1-\ga)}{1-\al(1-\ga)}>0\]
Then there is even enough space to fit in some small $c\ve>0$ for $\ve>0$ small enough, s.t.
\[(1-\al)(1+c')+\al\ga+c\ve<1.\]}
allows one to absorb the first term into the LHS of \eqref{eq:2.3.100.v} yielding
\begin{flalign}
\int\limits_{B_{2r}}{\eta^2(\xb)|\D^{h,s}\grad u|^2(x)\;dx}\le&\frac{C(\ga)}{r^2}\int\limits_{B_{2r}}{|\grad u(x)|^2\;dx},&
\end{flalign}
where the RHS is independent of $h.$ Finally, Nirenberg's lemma implies that $u\in W_{loc}^{2,2}(\Om,\R^2).$\vspace{1cm}

\begin{re}
1. A warning might be in order here. By standard theory one would expect that for a $2\ti2-$dimensional uniformly convex (that is $0<\ga<1$) integrand and $u\in W_{loc}^{2,2}(\Om,\R^2)$ that the integrability is automatically improved to $u\in W_{loc}^{2,p}(\Om,\R^2)$ for some $p>2$ by Meyers' theorem and even $u\in C_{loc}^\infty(\Om,\R^2)$ by Schauder theory. However, note that for our integrand Meyers' does not apply. Indeed, recall that Meyers' is applied in the following way: By $u\in W_{loc}^{2,2}(\Om,\R^2)$ one knows that $\grad u\in W_{loc}^{1,2}(\Om,\R^4)$ solves the `linearized' elliptic PDE 
\[
-\div(A(x)\grad^2 u)=0\;\mb{in}\;D'(\Om,\R^4),
\]
 where $A(\cd):=\grad_\xi^2W(\grad u(\cd))$ is the coefficients matrix. Note, that $\xi\mapsto W(\xi)$ is smooth hence the regularity of $A$ is determined by the regularity of $\grad u.$ Now one of the prerequisites of Meyers' theorem is that $A\in L_{loc}^\infty.$  But, $\xi\mapsto\grad_\xi^2W(\xi)$ grows quadratically in $\xi.$ Hence, at the present stage, we can only guarantee $A\in L^1_{loc},$ but not necessarily $A\in L_{loc}^\infty.$ Hence, we need to obtain higher integrability by foot, which  in fact can be done for a general $\ga>0$ with an additional assumption becoming mandatory if $\ga\ge1.$ This is the subject of the next section.\\

2. We do not expect that theorem \ref{thm:2.3.DM08} can be fully extended to the range when $\ga\ge1$ as the famous counterexample by Szekelyhidi \cite{Sz04} indicates. It might however be possible to show, that stationary points are of the class $W_{loc}^{2,q}(\Om,\R^2)$ for some $1\le q<2$ by altering the proof of theorem \ref{thm:2.3.DM08}. There is also hope that for global ($\omega-$)minimizers even more can be said.
\label{re:sec.3.1}
\end{re}

\section{Higher-order regularity and Schauder theory}
\label{sec:2.4}
The content of this section will be the proof of Theorem \ref{Thm:2.4.1}, the main argument of which can be found in §4.1.\! It is followed by some supporting content: §4.2.\! describes two technical lemmatas used in the proof. §4.3.\! is a version of the reverse Fatou's lemma and §4.3.\! discusses Hardy spaces, their dual spaces, and the compensated compactness results we rely on in the proof of theorem \ref{Thm:2.4.1}.

\begin{center}\textsc{4.1. Proof of theorem \ref{Thm:2.4.1}}\end{center}

Key idea of the proof is to establish a reverse Hölder inequality. The technicalities in the following might be overwhelming so it might be helpful to some if we outline the strategy here:
\begin{enumerate}
\item Differentiating the ELE.
\item Approximating the candidate $u$ by a smooth map $u_\ve$ and properties thereof.
\item Testing the ELE with the standard mollifier $\vp=\eta^2u_\ve,$ will split the ELE in a difficult part $(I)_\ve$ containing the highest order terms and a simpler one $(II)_\ve$ containing the mixed terms.
\item Convergence of $(I)_\ve$ is shown by making use of a version of the reverse Fatou's lemma.
\item Estimating $(I)_\ve$ by means of compensated compactness.
\item Convergence and upper bound on $(II)_\ve$ are shown via standard estimates.
 \item The proof is then concluded by limit taking, an application of Gehring's lemma, and Schauder theory.
\end{enumerate}

 We start by\\

\textbf{Step 1: differentiating the ELE.}
The information that $u$ is of class $W_{loc}^{2,2}$ allows one to take partial derivatives of $\grad_\xi W(\grad u):$ For $i,j,k\in\{1,2\}$ we get 
\begin{flalign}
\p_k(\grad_\xi W(\grad u))_{ij}=&\d_{ia}\d_{jb}\p_{k}\p_au_b+\rho'(d)\p_{\xi_{ab}}((\cof\xi)_{ij})\p_{k}\p_au_b&\nonumber\\
&+\rho''(d)(\cof\grad u)_{ij}(\cof\grad u)_{ab}\p_{k}\p_au_b&\nonumber\\
=&(\grad_{\xi}^2 W(\grad u))_{ijab}\p_{k}\p_au_b&
\label{eq: 2.HR.0}
\end{flalign}
where 
\[(\grad_{\xi}^2 W(\grad u))_{ijab}=\d_{ia}\d_{jb}+\rho'(d)\p_{\xi_{ab}}((\cof\xi)_{ij})+\rho''(d)(\cof\grad u)_{ij}(\cof\grad u)_{ab}.\]
This representation emphasises the quasilinear structure of this quantity on the LHS of \eqref{eq: 2.HR.0}. Furthermore, a priori by considering the $L^p-$spaces the first and second term on the RHS of \eqref{eq: 2.HR.0} are in $L^2,$ while the rightmost one and therefore the quantity on the LHS are in $L^q$ for all $1\le q<2$ but not necessarily in any better space. \\

Now for every $k\in \{1,2\}$ we can test against $\vp_k\in C_c^\infty(\Om,\R^{2})$ and get
\begin{equation}
\sum\limits_{i,j}\int\limits_\Om {\p_k(\grad_\xi W(\grad u))_{ij}\p_{i}\vp_{kj}\;dx}=0.
\label{eq: 2.HR.1}
\end{equation}
Instead of having a system of equations we can sum over all $k:$
\begin{equation}
\sum\limits_{i,j,k}\int\limits_\Om {\p_k(\grad_\xi W(\grad u))_{ij}\p_{i}\vp_{kj}\;dx}=0,
\label{eq: 2.HR.2}
\end{equation}
needs to be satisfied for arbitrary $\vp_k\in C_c^\infty(\Om,\R^{2}),$ $k\in \{1,2\}.$\footnote{Note, this is indeed the same. It's easy to see that \eqref{eq: 2.HR.2}  follows from \eqref{eq: 2.HR.1}. So assume that \eqref{eq: 2.HR.2} holds and for the sake of a contradiction, that \eqref{eq: 2.HR.1} is false. Wlog. there exists $\vp_1$ s.t.  $\int\limits_\Om {\p_1(\grad_\xi W(\grad u))_{ij}\p_{x_i}\vp_{1j}\;dx}\not=0.$ By \eqref{eq: 2.HR.2} $\int\limits_\Om {\p_2(\grad_\xi W(\grad u))_{ij}\p_{x_i}\vp_{2j}\;dx}\not=0$ for all $\vp_2\in C_c^\infty(\Om,\R^{2})$ in particular $\vp_2\equiv0,$ a contradiction.} As usual we suppress the sums from now on.\\

By introducing $\phi:=(\vp_1,\vp_2)\in C_c^\infty(\Om,\R^{2\ti2})$ we can rewrite this PDE in a closed form 
\[\int\limits_{\Om}\grad_{\xi}^2 W(\grad u)\grad^2 u\cd\grad\phi=0\;\mb{for all}\;\phi\in C_c^\infty(\Om,\R^{2\ti2}),\]
where the multiplication needs to be understood in the above sense. We want to emphasize again that the quantity $\grad_\xi^2W(\grad u(\cd))\grad^2 u(\cd)$ is in $L^q$ for all $1\le q<2$ but not necessarily in $L^2.$ Therefore, we can only test with functions in $W_0^{1,q'}$, where $2<q'<\infty$ is the dual Hölder exponent. This produces the first technical difficulty since we are initially not allowed to test with the standard mollifier $\vp_k=\eta^2\p_k u.$ This is dealt with by\\

\textbf{Step 2: approximating $u.$} Define $\tilde{u}:=u\chi_{B(x_0,2r)}$ and $u_\ve:=\eta_\ve*\tilde{u}.$ Recall, that $\supp\eta_\ve\ss B(0,\ve)$ and $\supp\tilde{u}\ss B(x_0,2r).$ Furthermore, recall, the following strong convergence properties of the mollification: 
\begin{align*}
&u_\ve\ra u \;\mb{in}\; W^{1,p} \;\mb{for all}\;1\le p<\infty,&\\
&\grad^2 u_\ve\ra \grad^2 u \;\mb{in}\; L^2,&\\
&\det\grad u_\ve\ra\det\grad u \;\mb{in}\; L^p\;\mb{for all}\;1\le p<\infty,&\\
&\rho^{(k)}(d_{\grad u_\ve})\ra\rho^{(k)}(d_{\grad u}) \;\mb{in}\; L^p\;\mb{for all}\;1\le p<\infty \;\mb{and any}\;k\in \N.&
\end{align*}
These properties can be found in \cite[Appendix C.4]{LE10} or are simple consequences thereof.\\

\textbf{Step 3: testing the equation.} For each $k\in \{1,2\}$ we are testing the system by
\[\vp_k=\eta^2\p_k u_\ve.\] 
Note, $u_\ve$ is smooth with compact support, hence $\vp_k\in C_c^\infty$ can be used as a test function in \eqref{eq: 2.HR.2} yielding
\begin{flalign*}
\int\limits_{B_{2r}}{\eta^2(x)\grad\p_k u\cd \grad\p_k u_\ve\;dx}=&-\int\limits_{B_{2r}}{\eta^2(x)\p_k(\rho'(d_{\grad u})\cof\grad u)\cd\grad\p_ku_\ve\;dx}&\\
&-2\int\limits_{B_{2r}}{\eta(x)\grad_\xi^2W(\grad u)\grad \p_ku\cd(\p_k u_\ve\ot\grad\eta)\;dx}&\\
=:&(I)_\ve+(II)_\ve.&
\end{flalign*}

\textbf{Step 4: convergence of $(I)_\ve$:}\\

\textbf{Weak global bound on $(I)_\ve$:}
Notice, again, that $\p_k(\rho'(d_{\grad u})\cof\grad u)\in L^q$ for $1\le q<2.$ Hence, for every $\ve>0$ is in $\p_k(\rho'(d_{\grad u})\cof\grad u)\cd\grad\p_ku_\ve\in L^q$ for $1\le q<2,$ since $\grad\p_ku_\ve\in L^\infty.$\\

As a first step we give a crude estimate to $(I)_\ve$ for small enough $\ve>0:$
Applying Lemma \ref{lem:2.HR.2}.(ii) yields 
\begin{flalign}
(I)_\ve=&-\int\limits_{B_{2r}}{\eta^2\p_k(\rho'(d_{\grad u})\cof\grad u)\cd\grad\p_ku_\ve\;dx}&\nonumber\\
\le&\ga\int\limits_{B_{2r}}{\eta^2\max\left\{
|\grad\p_s u|^2,\frac{1}{2}|\grad\p_s u|^2+\frac{1}{2}|\grad\p_s u_\ve|^2\right\}\;dx}&\nonumber\\
\le&\ga\int\limits_{B_{2r}\cap\Om_1}{\eta^2|\grad\p_s u|^2\;dx}+\frac{\ga}{2}\int\limits_{B_{2r}\cap\Om_2}{\eta^2(|\grad\p_s u|^2+|\grad\p_s u_\ve|^2)\;dx}&\nonumber\\
\le&\ga\int\limits_{B_{2r}}{|\grad\p_k u|^2\;dx},&
\label{eq:Ch1.4.101}
\end{flalign}
for small enough $\ve>0$ and where $\Om_1,\Om_2\ss\Om$ are disjoint sets s.t. $|\grad\p_s u|^2\ge|\grad\p_s u_\ve|^2$ a.e. in $\Om_1$ and $|\grad\p_s u|^2<|\grad\p_s u_\ve|^2$ a.e. in $\Om_2.$ \\
In the last step we used the following property of the mollification:\\
Let $V\ss\ss W\ss\ss U$ and $V,W,U$ open and bounded sets. Assume $f\in L_{loc}^p(U)$ and $1\le p<\infty.$ If $\ve>0$ small enough then for $\|f_\ve\|_{L^p(V)}\le\|f\|_{L^p(W)},$ see again \cite[Appendix C.4, p.631]{LE10}.\\
Indeed, we can set $U=\Om,$ $W=B_{2r}$ and we can find an open set $V$ s.t. $\supp\eta\ss\ss V\ss\ss B_{2r}.$ Moreover, choosing $p=2,$ $f:=\grad\p_k u$ yields the estimate.\\

Note that the above inequality is good enough in the case that $0<\ga<1,$ see Remark \ref{re:2.3.1} below. However, it is too crude if $\ga\ge1.$ Hence, we need to give a more refined argument. The strategy will be as follows: Instead of the crude estimate given in \eqref{eq:Ch1.4.101} we would like to consider the $\limsup\limits_{\ve\ra0}(I)_\ve$ and then get a more refined bound on this limit.\\

\textbf{Existence of $\limsup\limits_{\ve\ra0}(I)_\ve$:}
To guarantee the existence of $\limsup\limits_{\ve\ra0}(I)_\ve$ we want to apply a version of the reverse Fatou's lemma. The standard Fatou's lemma requires a pointwise majorant on the considered sequence. It is a subtle point, that $\ga |\grad\p_k u|^2$ turns out to be a global bound as shown in \eqref{eq:Ch1.4.101} but not a pointwise majorant for our sequence. Luckily there is a version of Fatou's, which only requires an pointwise individual integrable majorant for every member, as long as this sequence of majorants converges themself.
We are exactly in such a situation, to see this by lemma \ref{lem:2.HR.2}.(ii) we have
\begin{flalign}
-\eta^2\left(\p_k(\rho'(d_{\grad u})\cof\grad u)\cd\grad\p_su_\ve\right)(x)\le&\ga\max\left\{
|\grad\p_k u(x)|^2,\frac{1}{2}|\grad\p_k u(x)|^2+\frac{1}{2}|\grad\p_k u_\ve(x)|^2\right\}&
\label{eq:Ch1.4.1011}
\end{flalign}
for a.e.\! $x\in B_{2r}$ and for every $\ve>0.$
Define the sequence 
\begin{align*}g_\ve:=\ga\max\left\{
|\grad\p_k u|^2,\frac{1}{2}|\grad\p_k u|^2+\frac{1}{2}|\grad\p_k u_\ve|^2\right\}\in L^1(B_{2r},[0,\infty])\end{align*}
for every $\ve>0.$
Obviously, $g_\ve\ra \ga |\grad\p_k u|^2=:g$ converges strongly in $L^1(B_{2r},[0,\infty])$  and $g\in L^1(B_{2r},[0,\infty]).$ \\
Furthermore, define the sequence of integrable functions
\begin{align*}f_\ve:=-\eta^2\left(\p_k(\rho'(d_{\grad u})\cof\grad u)\cd\grad\p_ku_\ve\right)\in  L^1(B_{2r},[0,\infty])\end{align*} for every $\ve>0.$  Moreover, note  \begin{align*}f_\ve\ra f:=-\eta^2\left(\p_k(\rho'(d_{\grad u})\cof\grad u)\cd\grad\p_ku\right)\;\mb{pw. a.e. on}\;B_{2r}. \end{align*}
This convergence follows from the fact that $\grad\p_ku_\ve\ra\grad\p_ku$ converges pw. a.e. Note, that while the $f_\ve$'s are integrable, for the limit $f$ we can only guarantee that $f$ is measurable and $f^+$ integrable. Indeed, Lemma \ref{lem:2.HR.2}.(ii) guarantees an upper bound on $\int\limits_{B_{2r}}{f^+\;dx}$ given by 
\begin{align*}\int\limits_{B_{2r}}{f\;dx}\le\int\limits_{B_{2r}}{f^+\;dx}\le \ga\int\limits_{B_{2r}}{|\grad\p_k u|^2\;dx}<\infty.
\end{align*}
However, we do not have a lower bound on $\int\limits_{B_{2r}}{f\;dx}$ so it might be $-\infty.$ Luckily, our version of the reverse Fatou's lemma is such that it can still be applied.\\

Now after all those preparations we are finally in the position to apply our version of the reverse Fatou's lemma. By applying lemma \ref{lem:2.3.RF1} we get
\begin{align*}
\limsup\limits_{\ve\ra0}(I)_\ve=&\limsup\limits_{\ve\ra0}\int\limits_{B_{2r}}{f_\ve\;dx}&\\
\le&\int\limits_{B_{2r}}{f\;dx}&\\
=:&(I).&
\end{align*}

Now we want to bound $(I)$ from above. The first case we have to take care of is that $(I)$ might be $-\infty.$ In this case we can just estimate $(I)$ by $0.$\\

So assume from now on, that $(I)$ is finite. Making use of the product rule yields
\begin{align*}
(I)=&-\int\limits_{B_{2r}}{\eta^2\p_k(\rho'(d_{\grad u})\cof\grad u)\cd\grad\p_ku\;dx}&\\
=&-\int\limits_{B_{2r}}{\eta^2\rho'(d_{\grad u})\cof\grad u,_k\cd\grad u,_k\;dx}-\int\limits_{B_{2r}}{\eta^2\rho''(d_{\grad u})(\cof\grad u\cd\grad\p_ku)^2\;dx}.&
\end{align*}
Realising that the rightmost term is always non-positive implies
\begin{align*}
(I)\le&-\int\limits_{B_{2r}}{\eta^2\rho'(d_{\grad u})\cof\grad u,_k\cd\grad u,_k\;dx}.&
\end{align*}

So far we have seen, that the $\limsup\limits_{\ve\ra0}(I)_\ve$ is either $-\infty,$ which is harmless as argued, or it is finite and the explicit upper bound is given by 
\[
\limsup\limits_{\ve\ra0}(I)_\ve\le (I)\le(I.a):=-\int\limits_{B_{2r}}{\eta^2\rho'(d_{\grad u})\cof\grad u,_k\cd\grad u,_k\;dx}.\]
\vspace{0.5cm}

\textbf{Step 5: bounding $(I.a)$ via compensated compactness:}\\

\textbf{Approximating $(I.a):$} For reasons, becoming obvious in a second, we introduce, again, an approximation:
\begin{align}
(I.a)=&-\int\limits_{B_{2r}}{\eta^2\rho'(d_{\grad u})\cof\grad u_{\ve,k}\cd\grad u_{\ve,k}\;dx}&\nonumber\\
&+\int\limits_{B_{2r}}{\eta^2\rho'(d_{\grad u})\cof(\grad u_{\ve,k}-\grad u_{,k})\cd\grad u_{\ve,k}\;dx}&\nonumber\\
&+\int\limits_{B_{2r}}{\eta^2\rho'(d_{\grad u})\cof\grad u_{,k}\cd(\grad u_{\ve,k}-\grad u_{,k})\;dx}&\nonumber\\
=:&(I.a.a)_\ve+(I.a.b)_\ve+(I.a.c)_\ve.&
\end{align}
The terms $(I.a.b)_\ve$ and $(I.a.c)_\ve$ vanish when $\ve\ra0.$ Indeed by Hölder's we get
\begin{align*}|(I.a.c)_\ve|\le&\ga\int\limits_{B_{2r}}{|\cof\grad u_{,k}\cd(\grad u_{\ve,k}-\grad u_{,k})|\;dx}\\
\le&\ga\|\grad u_{,k}\|_{L^2}\|\grad u_{\ve,k}-\grad u_{,k}\|_{L^2}\ra0
\end{align*}
when $\ve\ra0,$ since $\grad^2u\in L^2$ and $\grad^2u_\ve\ra\grad^2u$ strongly in $L^2.$ One can argue similarly for $(I.a.b)_\ve.$\\

\textbf{Strategy for estimating $(I.a.a)_\ve$:} Recall now that the goal is to establish a reverse Hölder inequality. For this sake, we need to estimate $(I.a.a)_\ve$ in terms of the $L^p-$norms of the 2nd derivative for some $p<2$.\\
The integrand of $(I.a.a)_\ve$ contains the product $\cof\grad u_{\ve,k}\cd\grad u_{\ve,k}.$ Obviously, these quantities satisfy $\div\cof\grad u_{\ve,k}=0$ and $\mb{curl} \grad u_{\ve,k}=0$ in a distributional sense. This prompts the idea of making use of compensated compactness results, in particular, the famous div-curl lemma. Some difficulties arise here. First, both factors contain second derivatives of $u,$ which we would like to control in an $L^p$ norm with $p<2,$ this is done below. For this sake, the div-curl lemma \ref{lem:2.HR.1} needs to be applied with $s=t=p<2,$ which is possible although, the quantity $\cof\grad u_{\ve,k}\cd\grad u_{\ve,k}$ is only controlled in the weak space $\mathcal{H}^\al$ with $\al=\frac{p}{2}<1.$ Secondly, we need the dual space to such a space and the duality inequality, to bring the Hardy space into play in the first place.\\

\textbf{Extending the domain to $\R^2$:} For any map $v:B_{2r}\ra\R^2$ we will use $\ol{v}$ to denote the trivial extension, i.e. $\ol{v}:=v\chi_{B_{2r}},$ where $\chi$ is the characteristic function. Define $f:=\cof\grad u_{\ve,k}$ and  $g:=-\grad u_{\ve,k}.$ Then we can rewrite $(I.a.a)_\ve$ by
\begin{align}
(I.a.a)_\ve=&\int\limits_{\R^2}{\ol{\eta^2\rho'(d_{\grad {u}})}(\ol{f}\cd \ol{g})\;dx}.&\label{eq:Ch1.4.10221}
\end{align}
Note, that the extensions $\ol{f}$ and $\ol{g}$ are indeed smooth on the whole space. This is true because  $\grad u_{\ve,k}$ lives on a compact subset of $B_{2r}$ and disappears smoothly on the boundary of that compact set in such a way that it is smooth on the full set $\ol{B_{2r}}$ and remains zero up to the boundary of  $\ol{B_{2r}}.$ Moreover, they satisfy the cancelation conditions in $\ol{B_{2r}},$ hence, the extensions $\ol{f}$ and $\ol{g}$ are indeed smooth and satisfy the cancelation conditions $\div\ol{f}=0$ and $\mb{curl}\ol{f}=0$ on the full space.\\

As a next step we want to apply a Fefferman-Stein type duality inequality, however, for $\mathcal{H}^\al$ when $0<\al<1$ then the corresponding dual spaces 
are given by $\dot{\Lambda}_\beta/\{constants\}(\R^2),$ with $\beta=2(\frac{1}{\al}-1).$ We will need $H:=\ol{\eta^2\rho'(d_{\grad {u}})}\in \dot{\Lambda}_\beta(\R^2).$ 
As a first step $H$ must be continuous. For this first note, that $x\mt\rho'(d_{\grad {u(x)}})$ is continuous because of the assumption that  $x\mt d_{\grad u(x)}$ is Hölder-continuous.
Note, however, that $u$ might be non zero on the boundary of $B_{2r},$ hence $H$ might jump on the boundary. Obviously, $H$ is continuous in $\supp\eta$ and $H=0$ on $\R^2\sm \ol{B_{2r}}.$
However, $\eta$ disappears on $\ol{B_{2r}}\sm \supp\eta$ guaranteeing the continuity of $H.$ That $H$ is actually in $\dot{\Lambda}_\beta(\R^2)$ is shown below.\\

\textbf{Compensated compactness:} 
Applying now the duality inequality \eqref{eq:Ch1.4.103} to the RHS of \eqref{eq:Ch1.4.10221}, with a for now free $0<\al<1$ and the corresponding $\beta=2(\frac{1}{\al}-1),$ gives
\begin{align*}
(I.a.a)_\ve \le&C(\al) \|H\|_{\dot{\Lambda}_\beta(\R^2)}\|\ol{f}\cd \ol{g}\|_{\mathcal{H}^\al(\R^2)}.&
\end{align*}
Making use of the div-curl lemma \ref{lem:2.HR.1} with $\al=\frac{p}{2}\in(\frac{2}{3},1)$ and $s=t=p,$ and using \eqref{eq:2.4.HS.equiv} yields
\[\|\ol{f}\cd \ol{g}\|_{\mathcal{H}^\al(\R^2)}\le  C(p)\|\ol{f}\|_{L^p(\R^2)}\|\ol{g}\|_{L^p(\R^2)}\le C(p)\|\grad u_{\ve,k}\|_{L^p(B_{2r})}^2.\]

As a last step we need to show, that $H\in\dot{\Lambda}_\beta(\R^2)$ with $\beta=2(\frac{2}{p}-1)$ and bound the norm by a constant.\\

\textbf{Showing $H\in\dot{\Lambda}_\beta(\R^2)$ and estimating the norm}: 
If $x,y\in \R^2\sm\ol{B_{2r}}$ then $|H(x)-H(y)|=0.$ The case when $x\in \R^2\sm\ol{B_{2r}}$ and $y\in \ol{B_{2r}}$ or vice versa can be reduced to the case below. Indeed, letting $z$ be the point where the straight line connecting $x$ and $y$ hits $\partial B_{2r},$ then by the triangle inequality and $H(y)=H(z)=0$ we get 
\begin{align*}|H(x)-H(y)|\le|H(x)-H(z)|+|H(z)-H(y)|=|H(x)-H(z)|.\end{align*}
By the result below and the fact that $|x-z|\le|x-y|$ by choice we finally get
\begin{align*}|H(x)-H(y)|\le C|x-z|^\be\le C|x-y|^\be.\end{align*}

Let now $x,y\in \ol{B_{2r}}.$ Then we have
\begin{align}|H(x)-H(y)|\le \ga|\eta^2(x)-\eta^2(y)|+|\rho'(d_{\grad u})(x)-\rho'(d_{\grad u})(y)|.\label{eq:Ch1.4.10221.a}\end{align}
By $\rho\in C^\infty$ and the assumption $d_{\grad u}\in C_{loc}^{0,\be}$ for the rightmost term we get 
\begin{align*}|\rho'(d_{\grad u})(x)-\rho'(d_{\grad u})(y)|\le&\|\rho''\|_{L^\infty}|d_{\grad u}(x)-d_{\grad u}(y)|\le C \|\rho''\|_{L^\infty}|x-y|^\be.&\end{align*}

In order to establish an estimate on the first term of the RHS of \eqref{eq:Ch1.4.10221.a} we follow the strategy of Morrey's proof of the Dirichlet growth theorem.\\
Let $h:=\eta^2.$ Then by the triangle inequality for an arbitrary $w\in \R^2$ it holds
\[|h(x)-h(y)|\le|h(x)-h(w)|+|h(w)-h(y)|.\]
Denote the midpoint of $x$ and $y$ by $z,$ i.e. $z=\frac{x+y}{2}$ and the distance between $x$ and $z$ by $l:=\frac{|x-y|}{2}.$ Now we average the previous inequality over $B(z,l)$, i.e. \[|h(x)-h(y)|\le\fint\limits_{B(z,l)}{|h(x)-h(w)|\;dw}+\fint\limits_{B(z,l)}{|h(w)-h(y)|\;dw}.\]
It is enough to control one of them the other one is similar. For this sake, let $w_t:=w+t(x-w).$ By the mean value theorem we see
\begin{flalign*}
\fint\limits_{B(z,l)}{|h(x)-h(w)|\;dw}=&\fint\limits_{B(z,l)}{\left|\int\limits_0^1\grad h(w_t)dt\cd(x-w)\right|\;dw}&\\
\le&\fint\limits_{B(z,l)}\int\limits_0^1{|\grad h(w_t)|\;dt\;dw}|x-y|&\\
=&\frac{|x-y|}{|B(z,l)|}\int\limits_0^1\int\limits_{B(\bar{z},lt)}{|\grad h(v)|\;dv\;t^{-2}dt},&
\end{flalign*}
where we exchanged integrals and substituted $v=y+t(w-y)$ and $\bar{z}=y+t(z-y).$
Expressing the function explicitly, yields
\begin{align}
\int\limits_{B(\bar{z},lt)}{|\grad h(v)|\;dv}=&\int\limits_{B(\bar{z},lt)}{|2\eta(v)\p_k\eta(v)|\;dv}&\label{eq:CH1.4.105}
\end{align}
Applying Hölder's inequality for any $1<q<\infty$ we get 
\begin{align*}
\int\limits_{B(\bar{z},lt)}{|2\eta(v)\p_k\eta(v)|\;dv}\le &\frac{C}{r}\int\limits_{B(\bar{z},lt)\cap B(x_0,2r)}{\;dv}&\\
\le& \frac{C}{r}|B(\bar{z},lt)|^{\frac{1}{q}}|B(x_0,2r)|^{\frac{(q-1)}{q}}&\\
=&Cr^{1-\frac{2}{q}}(lt)^{\frac{2}{q}},&
\end{align*}
where we used $\eta\le 1,$ $|\grad\eta|\le \frac{c}{r},$ $|\rho'|\le\ga$ and $|B(\bar{z},lt)|=\pi(lt)^2.$
Finally,
\begin{flalign*}
|h(x)-h(y)|\le&\frac{|x-y|}{|B(z,l)|}Cr^{1-\frac{2}{q}}l^{\frac{2}{q}}\int\limits_0^1{t^{\frac{2}{q}-2}\;dt}&\\
\le&C(q)r^{1-\frac{2}{q}}|x-y|^{{\frac{2}{q}}-1}.
\end{flalign*}
 It is important to match $\beta$ exactly, since on the full space $\R^2$ there is no obvious relation between (homogenous) Lipschitz/Hölder-spaces with different exponents. Consequently, ${{\frac{2}{q}}-1}=\beta=2(\frac{2}{p}-1)$ needs to be satisfied. This choice yields $\frac{2}{q}=1+\beta>1$ guaranteeing that the integrals above can always be performed.\\

Hence, we showed that $H\in\dot{\Lambda}_{\beta}$ and the norm is bounded via
\[\|H\|_{\dot{\Lambda}_{\beta}}\le C(\ga,p)r^{-\beta}.\]

\textbf{Conclusion of Step 4+5:}
Finally, combining the above we obtain
\begin{align}
(I.a.a)_\ve\le C(\ga,p)r^{-\beta}\|\grad u_{\ve,k}\|_{L^p(B_{2r})}^2,
\label{eq:Ch1.4.1023}
\end{align}
yielding
\[
-\infty\le\limsup\limits_{\ve\ra0}(I)_\ve\le(I.a)\le\limsup\limits_{\ve\ra0}(I.a.a)_\ve\le C(\ga,p)r^{-\beta}\|\grad u_{,k}\|_{L^p(B_{2r})}^2.
\]
\vspace{0.5cm}

\textbf{Step 6: upper bound on $(II)_\ve$:}
Splitting up the integrals is fine since all the integrands are well behaved. Hence we get
\begin{flalign}
(II)_\ve=&-2\int\limits_{B_{2r}}{\eta(x)\grad_\xi^2W(\grad u)\grad \p_ku\cd(\p_k u_\ve\ot\grad\eta)\;dx}&\nonumber\\
=&-2\int\limits_{B_{2r}}{\eta(x)\grad \p_k u\cd(\p_k u_\ve\ot\grad\eta)\;dx}&\nonumber\\
&-2\int\limits_{B_{2r}}{\eta(x)\rho'(d_{\grad u})(\p_{\xi_{ab}}(\cof\xi)_{ij})\p_{k}\p_au_b\p_k (u_\ve)_i\p_j\eta\;dx}&\nonumber\\
&-2\int\limits_{B_{2r}}{\eta(x)\rho''(d)(\cof\grad u)_{ij}(\cof\grad u)_{ab}\p_{x_k}\p_au_b\p_k (u_\ve)_i\p_j\eta\;dx}&\nonumber\\
=:&(II.a)_\ve+(II.b)_\ve+(II.c)_\ve.&\label{eq:Ch1.4.99}
\end{flalign}

\textbf{Estimating $(II.a)_\ve+(II.b)_\ve$:} the first two terms $(II.a)_\ve$ and $(II.b)_\ve$ can be controlled by Hölder's inequality (for any $1<p<2$ and $p^{-1}+p'^{-1}=1$) and recalling the properties of $\eta$ to obtain
\begin{flalign}
(II.a)_\ve+(II.b)_\ve\le&C(\ga)\int\limits_{B_{2r}}{\eta|\grad\p_ku||\p_k u_\ve||\grad\eta|\;dx}&\nonumber\\
\le&C(\ga)\left(\int\limits_{B_{2r}}{|\eta|^p|\grad\p_ku|^p\;dx}\right)^{\frac{1}{p}}\left(\int\limits_{B_{2r}}{|\p_k u_\ve|^{p'}|\grad\eta|^{p'}\;dx}\right)^{\frac{1}{p'}}&\nonumber\\
\le&C(\ga)r^{-1}\left(\int\limits_{B_{2r}}{|\grad\p_ku|^p\;dx}\right)^{\frac{1}{p}}\left(\int\limits_{B_{2r}\sm B_r}{|\p_k u_\ve|^{p'}\;dx}\right)^{\frac{1}{p'}}.&\label{eq:Ch1.4.991}
\end{flalign}
Now the idea is to give an estimate of the form 
 \begin{align*}
\|\p_k u_\ve\|_{L^{p'}(B_{2r}\sm B_r)}\lesssim \|\grad\p_k u_\ve\|_{L^p(B_{2r})}.
 \end{align*}

For this sake, recall the Sobolev-Poincare inequality for balls\footnote{See, \cite[Theorem 1, §1.4.4.]{M11}.}: there is a constant $C(s^*)>0$ s.t. for any $v\in W_0^{1,s^*}$ it holds
 \begin{align}\|v\|_{L^{s}(B(x,R))}\le C(s^*)R^{1+\frac{2}{s}-\frac{2}{s^*}}\|\grad v\|_{L^{s^*}(B(x,R))},\label{eq:Ch1.4.992} \end{align}
where $s^*=\frac{2s}{2+s}.$ By applying the latter inequality with $R=2r,$ $s=p'$ and $s^*=(p')^*=\frac{2p'}{2+p'}$ we obtain 
 \begin{align*}
\|\p_k u_\ve\|_{L^{p'}(B_{2r})}\le& C(p)r^{1+\frac{2}{p'}-\frac{2}{(p')^*}}\|\grad\p_k u_\ve\|_{L^{(p')^*}(B_{2r})}.&
\end{align*}
Making use of Hölder's to go from $L^{(p')^*}$ to $L^{p}$ we get
 \begin{align*}
\|\p_k u_\ve\|_{L^{p'}(B_{2r})}\le C(p)r^{1+\frac{2}{p'}-\frac{2}{(p')^*}+\frac{2}{(p')^*}-\frac{2}{p}}\|\grad\p_k u_\ve\|_{L^p(B_{2r})}, \end{align*}
keeping in mind, that this can only be true, if the condition $(p')^*\le p$ is satisfied. Realise that the exponent in the latter inequality simplifies to 
 \begin{align*}\omega_1:=1+\frac{2}{p'}-\frac{2}{p}. \end{align*}
Hence, we have seen that
 \begin{align*}
\|\p_k u_\ve\|_{L^{p'}(B_{2r}\sm B_r)}\le\|\p_k u_\ve\|_{L^{p'}(B_{2r})}\le C(p)r^{\omega_1}\|\grad\p_k u_\ve\|_{L^p(B_{2r})}.
 \end{align*}
By combining the latter inequality with \eqref{eq:Ch1.4.991}  we obtain
 \begin{align}(II.a)_\ve+(II.b)_\ve\le C(\ga,p)r^{(\omega_1-1)}\|\grad\p_k u\|_{L^p(B_{2r})}\|\grad\p_k u_\ve\|_{L^p(B_{2r})}.\label{eq:Ch1.4.993} \end{align}
\vspace{0.25cm}

\textbf{Estimating $(II.c)_\ve$:} For the rightmost term in \eqref{eq:Ch1.4.99} we argue similarly as before, again by Hölder's inequality and the properties of $\eta$ one obtains
\begin{flalign*}
(II.c)_\ve=&-2\int\limits_{B_{2r}}{\eta(x)\rho''(d)(\cof\grad u)_{ij}(\cof\grad u)_{ab}\p_{k}\p_au_b\p_k (u_\ve)_i\p_j\eta\;dx}&\\
\le&C(\|\rho''\|_{L^\infty})\int\limits_{B_{2r}}{\eta(x)|\p_k u_\ve||\grad\p_ku||\grad u|^2|\grad\eta|\;dx}&\\
\le& C(\|\rho''\|_{L^\infty})r^{-1}\left(\int\limits_{B_{2r}}{|\grad \p_ku|^p\;dx}\right)^{\frac{1}{p}}\left(\int\limits_{B_{2r}\sm B_r}{|\p_k u_\ve|^{{p'}}|\grad u|^{{2p'}}\;dx}\right)^{\frac{1}{{p'}}}.&
\end{flalign*}
By applying Hölder's inequality with $\frac{2}{3}+\frac{1}{3}=1$ to the rightmost term of the latter line and by recalling the property $\|\p_k u_\ve\|_{L^{q}(B_{2r})}\le\|\p_k u\|_{L^{q}(B_{2r})}$ for any $1\le q<\infty$ and for small enough $\ve>0$ we get 
\begin{flalign*}
\left(\int\limits_{B_{2r}\sm B_r}{|\p_k u_\ve|^{{p'}}|\grad u|^{{2p'}}\;dx}\right)^{\frac{1}{{p'}}}\le& \left(\int\limits_{B_{2r}\sm B_r}{|\grad u|^{{3p'}}\;dx}\right)^{\frac{2}{{3p'}}}\left(\int\limits_{B_{2r}\sm B_r}{|\p_k u_\ve|^{{3p'}}\;dx}\right)^{\frac{1}{{3p'}}}&\\
\le&\left(\int\limits_{B_{2r}\sm B_r}{|\grad u|^{{3p'}}\;dx}\right)^{\frac{1}{{p'}}},&
\end{flalign*}
if $\ve>0$ is sufficiently small.\\

We can proceed like before. Applying again \eqref{eq:Ch1.4.992} yields 
\begin{align*}
\|\grad u\|_{L^{3p'}(B_{2r})}\le& C(p)r^{1+\frac{2}{3p'}-\frac{2}{(3p')^*}}\|\grad^2 u\|_{L^{(3p')^*}(B_{2r})}.&
\end{align*}
Assuming $(3p')^*\le p$ and applying Hölder's to go from $L^{(3p')^*}$ to $L^{p}$ we get
 \begin{align*}
\|\grad u\|_{L^{3p'}(B_{2r})}\le C(p)r^{1+\frac{2}{3p'}-\frac{2}{(3p')^*}+\frac{2}{(3p')^*}-\frac{2}{p}}\|\grad^2 u\|_{L^p(B_{2r})}. \end{align*}
 The exponent in the latter inequality simplifies to 
 \begin{align*}\omega_2:=1+\frac{2}{3p'}-\frac{2}{p}. \end{align*}
This shows
 \begin{align*}
\|\grad u\|_{L^{3p'}(B_{2r}\sm B_r)}^3\le\|\grad u\|_{L^{3p'}(B_{2r})}^3\le C(p)r^{3\omega_2}\|\grad^2 u\|_{L^p(B_{2r})}^3.
\end{align*}
Then
 \begin{equation}
 (II.c)_\ve\le C(\|\rho''\|_{L^\infty},p)r^{\left(3\omega_2-1\right)}\|\grad\p_k u\|_{L^p(B_{2r})}\|\grad^2 u\|_{L^p(B_{2r})}^3. 
 \label{eq:Ch1.4.994} 
 \end{equation}
\vspace{0.25cm}

\textbf{Estimating $(II)_\ve$:} 
By combining \eqref{eq:Ch1.4.993} and \eqref{eq:Ch1.4.994} we have
\[(II)_\ve\le C(\ga,\|\rho''\|_{L^\infty},p)\|\grad\p_k u\|_{L^p(B_{2r})}\left[r^{(\omega_1-1)}\|\grad\p_k u_\ve\|_{L^p(B_{2r})}+r^{(3\omega_2-1)}\|\grad^2 u\|_{L^p(B_{2r})}^3\right],\]
for any $1<p<2$ s.t. $(p')^*\le p$ and $(3p')^*\le p,$ where the latter condition is stricter and reduces the range to $p\in[\frac{8}{5},2).$ 
\vspace{0.5cm}

\textbf{Step 7: conclusion:}\\

Taking the limsup $\ve\ra0$ on both sides, and collecting all expressions from above yields

\[
\int\limits_{B_{r}}{|\grad^2 u|^2\;dx}\le C(\ga,\|\rho''\|_{L^\infty},p)[(r^{-\be}+ r^{(\omega_1-1)})\|\grad^2u\|_{L^p}^2+r^{(3\omega_2-1)} \|\grad^2 u\|_{L^{p}}^4]
\]
Recalling $\beta=2(\frac{2}{p}-1),$  $\omega_1-1=-\be,$ and dividing the latter by $|B_{2r}|$ yields
\begin{equation}
\fint\limits_{B_{r}}{|\grad^2 u|^2\;dx}\le C(\ga,\|\rho''\|_{L^\infty},p)(1+r^{(3\omega_2-1+\beta)} \|\grad^2 u\|_{L^{p}(B_{2r})}^2)\left(\fint\limits_{B_{2r}}{|\grad^2 u|^p\;dx}\right)^{\frac{2}{p}}.
\label{eq:Ch1.4.99.99}
\end{equation}
Now, notice firstly $3\omega_2-1+\be=\frac{2}{p'}-\frac{2}{p}$ and secondly by Hölder's inequality
 \begin{align*}
r^{\left(\frac{2}{p'}-\frac{2}{p}\right)}\|\grad^2 u\|_{L^{p}(B_{2r})}^2\le C(p)r^{\left(\frac{2}{p'}-\frac{2}{p}+\be\right)}\|\grad^2 u\|_{L^2(B_{2r})}^2.
\end{align*}
Then the exponent $\frac{2}{p'}-\frac{2}{p}+\be=\frac{2}{p'}+\frac{2}{p}-2=0$ vanishes for any $p\in[\frac{8}{5},2).$ Furthermore, by Lebesgue's absolute continuity theorem we can choose $r\le1$ so small s.t. $\|\grad^2 u\|_{L^2(B_{2r})}^2\le1$ and \eqref{eq:Ch1.4.99.99} becomes\footnote{A similar reasoning can be found in \cite[p.142]{MG83}.}
\begin{equation}
\fint\limits_{B_{r}}{|\grad^2 u|^2\;dx}\le C(\ga,\|\rho''\|_{L^\infty},p)\left(\fint\limits_{B_{2r}}{|\grad^2 u|^p\;dx}\right)^{q},
\label{eq:Ch1.4.100}
\end{equation}

where $q=\frac{2}{p}.$ By choosing $f:=|\grad^2 u|^p,$ Gehring's lemma, see \cite[Ch.5, Prop.1.1]{MG83} guarantees the existence of some $\d=\d(\ga,\|\rho''\|_{L^\infty},p)>0$ and $\tilde{C}=\tilde{C}(\ga,\|\rho''\|_{L^\infty},p)>0$ s.t.  
\[
\left(\fint\limits_{B_{r}}{|\grad^2 u|^{2+\d}\;dx}\right)^{\frac{1}{2+\d}}\le \tilde{C}\left(\fint\limits_{B_{2r}}{|\grad^2 u|^2\;dx}\right)^{\frac{p}{2}}
\]
and in particular $u\in W_{loc}^{2,2+\d}(\Om,\R^2)$ and by Sobolev embedding 
$u\in C_{loc}^{1,\al},$ for some $\al>0.$ Smoothness is then obtained via Schauder theory, consult \cite[Chapter 3, Theorem 3.(1-3)]{MG83}. \vspace{0.5cm}
\begin{re}
Note, that for the uniformly convex case $(0<\ga<1)$ the argument can be simplified and the assumption $x\mt d_{\grad u(x)}\in C_{loc}^{0,\be}(\Om)$ can be dropped, in particular, the difficult limit taking and compensated compactness arguments can be avoided. The estimate found in \eqref{eq:Ch1.4.101} is good enough if $(0<\ga<1).$ Indeed, combined with the bound on $(II)_\ve$ from above one gets the slightly weaker estimate
\[
\fint\limits_{B_{r}}{|\grad^2 u|^2\;dx}\le C(\ga,\|\rho''\|_{L^\infty},p)\left(\fint\limits_{B_{2r}}{|\grad^2 u|^p\;dx}\right)^{q}+\ga\fint\limits_{B_{2r}}{|\grad^2 u|^2\;dx},
\]
replacing \eqref{eq:Ch1.4.100}.
 However, Gehring's lemma still applies and one can conclude as above.
\label{re:2.3.1}
\end{re}\vspace{0.5cm}

\begin{center}\textsc{4.2. Technical lemmatas}\end{center}

In this paragraph we show two technical lemmatas, which were used in section 3 and §4.3. Initially, we establish important upper and lower bounds on quantities arising from differencing the non-linear quantity $\rho'(d_u)\cof\grad u.$ These bounds are important for the result obtained in section 3 however it is also a prerequisite for the result we obtain next.

\begin{lem} Let $h\in\R^+,$ $s\in\{0,1\},$ $u\in W^{1,2}(\Om,\R^2)$ and $\rho$ be as before. Then\\
(i) the lower bound
\begin{equation}
\Delta^{h,s}(\rho'(d)\cof\grad u)(x)\cdot\Delta^{h,s}\grad u(x)\ge-\ga|\Delta^{h,s}\grad u(x)|^2
\label{eq:2.6}
\end{equation}
and\\
(ii) the upper bound
\begin{equation}
|\Delta^{h,s}(\rho'(d)\cof\grad u)(x)|\le\ga|\Delta^{h,s}\grad u(x)|+\frac{2\ga}{h}|\grad u(x)|
\label{eq:2.16}
\end{equation}
hold a.e. in $\Om_{h}:=\{x\in\Om:\dist(x,\p\Om)>h\}.$
\label{HR:LUBs}
\end{lem}
\begin{proof}
(i) Let $x\in \Om_h,$ fix an $s\in \{0,1\}$ and define $\bar{x}:=x+he_s.$ Recalling the product rule for difference quotients, for $f:\Om_h\rightarrow \R$ and $g:\Om\rightarrow \R$ it holds that
\begin{equation}
\D^{h,s}(fg)(x)=\D^{h,s}(f)(x)g(x)+f(\bar{x})\D^{h,s}g(x).
\label{eq:2.3.100PRFD}
\end{equation}
We apply this product rule to obtain 
\begin{eqnarray}
\Delta^{h,s}(\rho'(d)\cof\grad u)(x)\cdot\Delta^{h,s}\grad u(x)&=&\Delta^{h,s}(\rho'(d))(x)(\cof\grad u)(x)\cdot\Delta^{h,s}\grad u(x)\nonumber\\
&&+\rho'(d)(\bar{x})\Delta^{h,s}(\cof\grad u)(x)\cdot\Delta^{h,s}\grad u(x)\nonumber\\
&=:&I+II.
\label{eq:MP,PR2}
\end{eqnarray}
By linearity of the gradient and the cofactor in two dimensions it holds that $\Delta^{h,s}\grad u(x)=\grad\Delta^{h,s} u(x)$ and $\Delta^{h,s}(\cof\grad u)(x)=\cof\Delta^{h,s}\grad u(x)$ for all $x\in \Om_h.$ Let us now introduce the following notation $A:=\grad u(\bar{x}),$ $B:=\grad u(x).$ Note that for every matrix $N\in \R^{2\ti2},$ the determinant satisfies the following scaling behaviour $d_{N/h}=h^{-2}d_N$ for $h\in \R\setminus\{0\}.$ Using these properties and definitions to rewrite the second term of \eqref{eq:MP,PR2} as 
\begin{eqnarray}
II&=&\rho'(d)(\bar{x})\Delta^{h,s}(\cof\grad u)(x)\cdot\Delta^{h,s}\grad u(x)\nonumber\\
&=&\rho'(d)(\bar{x})\cof\Delta^{h,s}\grad u(x)\cdot\Delta^{h,s}\grad u(x)\nonumber\\
&=&h^{-2}\rho'(d_A)\cof (A-B) \cdot(A-B)\nonumber\\
&=&2h^{-2}\rho'(d_A)d_{A-B}.
\label{eq:2.7}
\end{eqnarray}
The first term of (\ref{eq:MP,PR2}) takes the form
\begin{eqnarray}
I=h^{-1}\Delta^{h,s}(\rho'(d))(x)\cof B\cdot(A-B).
\label{eq:2.8}
\end{eqnarray}
For two matrices $N,N'\in\R^{2\times2}$ the determinant of $N+N'$ can be expanded in the following way
\begin{eqnarray}
d_{N+N'}=d_{N'}+d_{N}+\cof N\cdot N'.
\label{eq:2.9}
\end{eqnarray}
Now choose $N'=A-B$ and $N=B$ and we get 
\begin{eqnarray}
d_{A}=d_{A-B}+d_{B}+\cof B\cdot (A-B)
\label{eq:2.10}
\end{eqnarray}
or 
\begin{eqnarray}
\cof B\cdot (A-B)=d_{A}-d_{B}-d_{A-B}=h\Delta^{h,s}d(x)-d_{A-B}.
\label{eq:2.11}
\end{eqnarray}
Therefore,
\begin{eqnarray}
I=\Delta^{h,s}(\rho'(d))(x)\Delta^{h,s}d(x)-h^{-1}\Delta^{h,s}(\rho'(d))(x)d_{A-B}.
\label{eq:2.12}
\end{eqnarray}
The whole expression can now be written as
\begin{eqnarray}
I+II&=&\Delta^{h,s}(\rho'(d))(x)\Delta^{h,s}d(x)\nonumber\\
&&+h^{-2}(2\rho'(d_A)-(\rho'(d_A)-\rho'(d_B))d_{A-B}\nonumber\\
&=&\Delta^{h,s}(\rho'(d))(x)\Delta^{h,s}d(x)+h^{-2}(\rho'(d_A)+\rho'(d_B))d_{A-B}\nonumber\\
&=:&III+IV.
\label{eq:MP,PR3}
\end{eqnarray}
The latter term can be bounded from below through
\begin{eqnarray}
IV&\ge& -h^{-2}|\rho'(d_A)+\rho'(d_B)||d_{A-B}|\nonumber\\
&\ge&-h^{-2}\ga|A-B|^2,
\label{eq:2.13}
\end{eqnarray}
where we used $\rho'\le\ga$ and $|d_N|\le\frac{1}{2}|N|^2.$ Furthermore, realize that $III\ge0$  is always true, since $\rho$ is convex and $\rho'$ satisfies a monotonicity inequality. Indeed, consider first $\Delta^{h,s}d(x)\ge0$ then $d_A\ge d_B$ and by monotonicity of $\rho',$ $\Delta^{h,s}\rho'(d)(x)\ge0.$ Hence, $III\ge0.$ Secondly, let $\Delta^{h,s}d(x)<0$ then $d_A< d_B$ implying $\rho'(d_A)\le\rho'(d_B)$ and $\Delta^{h,s}\rho'(d)(x)\le0$ and again $III\ge0.$

This implies
\begin{equation}
\Delta^{h,s}(\rho'(d)\cof\grad u)\cdot\Delta^{h,s}\grad u\ge-\ga|\Delta^{h,s}\grad u|^2.
\label{eq:2.14}
\end{equation}
\\

(ii) The upper bound follows again from the product rule \eqref{eq:2.3.100PRFD} and the triangle inequality
\begin{eqnarray}
|\Delta^{h,s}(\rho'(d)\cof\grad u)(x)|&=&|\rho'(d)(\bar{x})\Delta^{h,s}(\cof\grad u)(x)+\Delta^{h,s}(\rho'(d))(x)(\cof\grad u)(x)|\nonumber\\
&\le&\ga|\Delta^{h,s}\grad u(x)|+\frac{2\ga}{h}|\grad u(x)|.
\end{eqnarray}\end{proof}
\vspace{1cm}

Now we obtain a lower bound on the quantity $\p_s(\rho'(d_{\grad v})\cof\grad v)\cd\grad \p_sv_\ve$ we used in §4.1. This follows from the lower bound we showed for the `differenced' version in Lemma \ref{HR:LUBs}.

\begin{lem}\label{lem:2.HR.2} Let $\Om\ss\R^2$ be open and bounded and $v\in W^{2,2}(\Om).$\\

(i). Then for $s=1,2$ it holds
\begin{equation}\p_s(\rho'(d_{\grad v})\cof\grad v)\cd\grad \p_sv\ge-\ga|\grad\p_s v|^2\end{equation}
a.e. in $\Om.$\\
(ii). For every $s=1,2$ and $\ve>0$ small enough it holds 
\begin{equation}\p_s(\rho'(d_{\grad v})\cof\grad v)\cd\grad \p_sv_\ve\ge-\ga\max\{|\grad\p_s v|^2,\frac{1}{2}|\grad\p_s v|^2+\frac{1}{2}|\grad\p_s v_\ve|^2\}.\end{equation}
a.e. in $\Om.$
\end{lem}
\begin{proof}
(i). By Nirenberg's lemma \cite[p.45]{MG83} we know that $\D^{h,s}\grad v\ra\grad\p_sv$ in $L^2(\Om,\R^4),$ in particular,
 $\D^{h,s}\grad v\ra\grad\p_sv$ a.e. in $\Om.$ By the standard/discrete product rule it is straightforward  to see that
\[\p_s(\rho'(d_{\grad v})\cof\grad v)=\lim\limits_{h\ra0}\D^{h,s}(\rho'(d_{\grad v})\cof\grad v)\]
a.e. in $\Om.$ By Lemma \ref{HR:LUBs} we get
\[\p_s(\rho'(d_{\grad v})\cof\grad v)=\lim\limits_{h\ra0}\D^{h,s}(\rho'(d_{\grad v})\cof\grad v)\ge-\ga\lim\limits_{h\ra0}|\D^{h,s}\grad v|^2=-\ga|\grad\p_s v|^2\]
a.e. in $\Om.$\\ 
(ii). There are two possible cases. Assume first that the inequality in (i) is strict, i.e.
\[\p_s(\rho'(d_{\grad v}(x))\cof\grad v(x))\cd\grad \p_sv(x)>-\ga|\grad\p_s v(x)|^2,\]
for some $x\in\Om.$
Since 
\[\p_s(\rho'(d_{\grad v}(x))\cof\grad v(x))\cd\grad \p_sv_\ve(x)\ra\p_s(\rho'(d_{\grad v}(x))\cof\grad v(x))\cd\grad \p_sv(x)\]
Then the same holds true if $\ve$ is small enough:\footnote{Indeed, we have a converging sequence $f_\ve\ra f$ if $\ve\ra 0$ with $f>g.$ Then $f_\ve>g$ for small enough $\ve.$ Assume not, then there exists a sequence $(\ve_n)_{n\in\N}>0$ with $\ve_n\ra 0$ if $n\ra\infty$ s.t. $f_{\ve_n}\le g$ this contradicts immediately the convergence $f_\ve\ra f$.}
\[\p_s(\rho'(d_{\grad v}(x))\cof\grad v(x))\cd\grad \p_sv_\ve(x)>-\ga|\grad\p_s v(x)|^2.\]

Assume instead that for some $x\in \Om$ equality holds:
\[\p_s(\rho'(d_{\grad v}(x))\cof\grad v(x))\cd\grad \p_sv(x)=-\ga|\grad\p_s v(x)|^2.\]
Then since, the above inequality is a Cauchy-Schwarz inequality, equality only holds if 
\[\p_s(\rho'(d_{\grad v}(x))\cof\grad v(x))=-\ga\grad\p_s v(x).\]
But then for every $\ve>0$ we get 
\begin{align*}\p_s(\rho'(d_{\grad v}(x))\cof\grad v(x))\cd\grad \p_sv_\ve(x)=&-\ga\grad\p_s v(x)\cd\grad \p_sv_\ve(x)&\\
\ge&-\frac{\ga}{2}(|\grad\p_s v(x)|^2+|\grad\p_s v_\ve(x)|^2).&
\end{align*}
Putting both cases together yields the claim.\end{proof}\vspace{0.5cm}

\begin{center}\textsc{4.3. A measure-theoretic result}\end{center}

The next statement is a generalisation of the the reverse Fatou's lemma. Instead of a pointwise integrable upper limit dominating every member of the sequence, it will be enough to bound each member individually from above, as long as this sequence of upper limits converges by itself. Our proof is inspired by \cite[Thm 2.1]{FKZ16}.

\begin{lem} [A version of the reverse Fatou's lemma]
Let $(X,\Sigma,\mu)$ be a measure space, with $\mu(X)<\infty.$ Let $g_n, g\in L^1(X,[-\infty,\infty],\mu)$ for any $n\in \N$ be a sequence, s.t.
\begin{align*}g_n \ra g \;\mb{strongly in} \; L^1(X,[-\infty,\infty],\mu).\end{align*}
Suppose further that there is another sequence $(f_n)_{n\in \N}\ss L^1(X,[-\infty,\infty],\mu)$ s.t. 
\begin{align*}f_n\le g_n \;\mb{a.e. for every} \;n\in\N\end{align*}
and s.t. there exists a measurable functions $f:X\ra[-\infty,\infty]$ 
s.t. \begin{align*}f_n\ra f\mb{ptw. a.e.}\end{align*} 
and either (a) $\int\limits_{X}{f\;d\mu}=+\infty$ or (b) $f^+$ is integrable.  \\

Then it holds
\begin{align}\limsup\limits_{n\ra\infty}\int\limits_{X}{f_n\;d\mu}\le \int\limits_{X}{f\;d\mu}.\end{align}
\label{lem:2.3.RF1}
\end{lem}
\begin{proof}
(a) $\int\limits_{X}{f\;d\mu}=+\infty.$ Trivial, since $\limsup\limits_{n\ra\infty}\int\limits_{X}{f_n\;d\mu}\le+\infty$ is always true.\\
(b) Assume $f^+$ to be integrable. \\
Note that we prove the more classical $\liminf-$version. This is done be making the necessary changes:  $F:=-f,F_n:=-f_n, G_n=-g_n, G=-g.$ Then all the statements above are still valid except $F_n\ge G_n$ and $F^-$ is now integrable. Then it is enough to show
\begin{align}\liminf\limits_{n\ra\infty}\int\limits_{X}{F_n\;d\mu}\ge \int\limits_{X}{F\;d\mu}.\label{eq:2.3.RF1.1}\end{align}

We start by considering the following decomposition
\begin{align*}\int\limits_{X}{F_n\;d\mu}=\int\limits_{X\cap\{F_n<-K\}}{F_n\;d\mu}+\int\limits_{X\cap\{F_n\ge-K\}}{F_n\;d\mu}=:A_{n,K}+B_{n,K},\end{align*}
which holds for any $n\in\N$ and any $K\in \R.$
The proof is divided into two steps. First one needs to show that $A_{n,K}\ra0$ for large enough $n$ and  $K.$ The second one is to show that $B_{n,K}$ can be controlled from below by $\int\limits_{X}{F\;d\mu}$ for large enough $n$ and  $K.$ Together with Step 1 that implies \eqref{eq:2.3.RF1.1}. \\

Step 1: $\liminf\limits_{K\ra+\infty}\liminf\limits_{n\ra+\infty} A_{n,K}\ge0.$\\

Since, $F_n\ge G_n$ we have
\begin{align*}A_{n,K}\ge& \int\limits_{X}{G_n\chi_{\{F_n<-K\}} \;d\mu}&\\
\ge& \int\limits_{X}{G_n\chi_{\{G_n<-K\}} \;d\mu}&
\end{align*}
for any $n\in\N$ and $K\in \R$. 
We can now split up the latter integral as follows
\begin{align}\int\limits_{X}{G_n\chi_{\{G_n<-K\}} \;d\mu}=\int\limits_{X}{(G_n-G)\chi_{\{G_n<-K\}} \;d\mu}+\int\limits_{X}{G\chi_{\{G_n<-K\}} \;d\mu}.\label{eq:2.3.RF1.2}\end{align}
Then the first term on the RHS converges to $0$ when $n\ra\infty,$ since  
\begin{align*}\int\limits_{X}{|(G_n-G)\chi_{\{G_n<-K\}}| \;d\mu}\le\int\limits_{X}{|G_n-G|\;d\mu}\ra0.\end{align*}
We can treat the rightmost term of \eqref{eq:2.3.RF1.2} by splitting up the integral as follows
\begin{align}\int\limits_{X}{G\chi_{\{G_n<-K\}} \;d\mu}=\int\limits_{X}{G\chi_{\{G_n<-K,|G_n-G|\le1\}} \;d\mu}+\int\limits_{X}{G\chi_{\{G_n<-K,|G_n-G|>1\}} \;d\mu}.\label{eq:2.3.RF1.3}\end{align}
Since $G\in L^1$ the rightmost term disappears when $n\ra\infty.$ Indeed, we have \\
\begin{align*}\mu(\{G_n<-K,|G_n-G|>1\})\le\mu(\{|G_n-G|>1\})\le\int\limits_{\{|G_n-G|>1\}}{|G_n-G| \;d\mu}\ra0\end{align*} when $n\ra\infty.$ Since $G\in L^1$ we have
\begin{align*}\int\limits_{X}{G\chi_{\{G_n<-K,|G_n-G|>1\}} \;d\mu}\ra0\end{align*}
when $n\ra\infty.$ \\
We can treat the first term of \eqref{eq:2.3.RF1.3} as follows. Realize first that
\begin{align*}\{G_n<-K,|G_n-G|\le1\}\ss\{G<-K+1\}\end{align*}
and hence by the Dominated Convergence Theorem we have
\begin{align*}G\chi_{\{G_n<-K,|G_n-G|\le1\}}\ra0 \;\mb{strongly in}\; L^1(X,[-\infty,\infty],\mu).\end{align*}
Together this shows
\begin{align*}\liminf\limits_{K\ra+\infty}\liminf\limits_{n\ra+\infty} A_{n,K}\ge0\end{align*}
completing Step 1.\\
 
Step 2:  $\liminf\limits_{K\ra+\infty}\liminf\limits_{n\ra+\infty} B_{n,K}\ge\int\limits_{X}{F\;d\mu}.$ \\
Since, $\mu(X)<\infty$ by Egorov's Theorem we can find for any $\ve>0$ we can find $X^\ve\ss X$ s.t. $\mu(X\sm X^\ve)<\ve$ and for any $\eta>0$ there is $N_\eta>0$ s.t. for any $n\ge N_\eta$ we have  $|F_n-F|<\eta$ on $X^\ve.$ With this at hand we can split up the integrals again as follows
\begin{align*}B_{n,K}=\int\limits_{X^\ve\cap\{F_n\ge-K\}}{F_n\;d\mu}+\int\limits_{(X\sm X^\ve)\cap\{F_n\ge-K\}}{F_n\;d\mu}.\end{align*}
For the rightmost term we find the simple estimate
\begin{align*}\int\limits_{(X\sm X^\ve)\cap\{F_n\ge-K\}}{F_n\;d\mu}\ge-K\mu(X\sm X^\ve)\ge-K\ve\end{align*} for any $K\in\R.$ Then since we have $|F_n-F|<\eta$ on $X^\ve$ this implies
\begin{align*}B_{n,K}\ge\int\limits_{X^\ve\cap\{F_n\ge-K\}}{F_n\;d\mu}-K\ve
\ge\int\limits_{X^\ve\cap\{F_n\ge-K\}}{(F-\eta)\;d\mu}-K\ve
\end{align*}
for any $n\ge N_\eta.$ Now we take $\ve_K:=K^{-2}$ and $\eta:=K^{-1}$ and $n\ge N(K)$.
Define $\phi_{n,K}:=F\chi_{\{X^{(\ve_K)}\cap\{F_n\ge-K\}\}}.$
Elementary set theoretic considerations yield for arbitrary sets $U,V,W$  with $V\ss U$ it holds
\begin{align*}U\sm(V\cap W)=(U\sm V)\cup (V\sm W)\ss(U\sm V)\cup (U\sm W).\end{align*}
Then for any $n\in \N$ we have
\begin{align*}\mu(X\sm\{X^{(\ve_K)}\cap\{F_n\ge-K\}\})\le\mu(X\sm X^{(\ve_K)})+\mu(X\sm\{F_n\ge-K\})\ra0\end{align*}
when $K\ra+\infty.$\\

Since  $F_n\ge G_n$ for any $n\in \N$ the same holds for the limit $F\ge G$ and hence we have $\phi_{n,K}=F\chi_{\{X^{(\ve_K)}\cap\{F_n\ge-K\}\}}\ge G$ for any $n\in \N$ and $K\in\R.$
Then for any $n\in \N,$ $(\phi_{n,K})_{K\in \R}$ is a sequence of measurable functions, which is bounded from below and we can apply the standard version of Fatou's Lemma to obtain\footnote{see \cite[Cor 5.34]{WZ15}.}
\begin{align*}\liminf\limits_{K\ra+\infty}\int\limits_{X}{\phi_{n,K}\;d\mu}\ge\int\limits_{X}{F\;d\mu}\;\mb{for any}\; n\in \N.\end{align*}
By putting everything together we get
\begin{align*}\liminf\limits_{K\ra+\infty}\liminf\limits_{n\ra+\infty}B_{n,K}\ge \int\limits_{X}{F\;d\mu},
\end{align*}
completing the proof.\end{proof}\vspace{0.5cm}

\begin{center}\textsc{4.4. Hardy spaces and duality}\end{center}
Initially, we recall the definition of Hardy spaces, where we follow \cite{U01}.
\begin{de}[Hardy spaces $\mathcal{H}^r$]
Let $S(\R^n)$ be the Schwartz class and $S'(\R^n)$ be its dual. For every $f\in S'(\R^n),$ and for  
$\Phi\in S(\R^n)$ s.t. $\int{\Phi \;dx}\not=0,$ a positive real parameter $\zeta>0$ and $x\in\R^n$ we define
 \[N_{\Phi,\zeta}f:=\sup\{f*(\Phi)_t(y):|x-y|<\zeta t\}.\]
Then for $r>0$ we can define the Hardy space by
\[\mathcal{H}_{\Phi,\zeta}^r(\R^n):=\{f\in S'(\R^n)|N_{\Phi,\zeta}f\in L^r(\R^n)\}.\]
Next, define the `norm'  by
\[\|f\|_{\mathcal{H}^r(\R^n)}:=\|N_{\Phi,\zeta}f\|_{L^r(\R^n)}.\]
Obviously, this quantity only defines a norm if $r\ge1$ and $\mathcal{H}^r$ is a Banach space. However, when $0<r<1,$ $\mathcal{H}^r(\R^n)$ is just a metric space. Hence, the notion of a dual space does not make any sense, in the later case.\\
It is important to note that these spaces are independent of $\Phi>0$ and $\zeta>0,$ i.e. for any $(\Phi,\zeta)$ and $(\Phi',\zeta')$ it holds 
\[\mathcal{H}_{\Phi,\zeta}^r(\R^n)=\mathcal{H}_{\Phi',\zeta'}^r(\R^n),\]
with equivalent `norms',i.e. for any choices of $(\Phi,\zeta)$ and $(\Phi',\zeta')$  there are constants $c=c(\Phi,\zeta)>0$ and $C=C(\Phi,\zeta)>0$ s.t. it holds
\[c\|N_{\Phi',\zeta'}f\|_{L^r(\R^n)}\le\|N_{\Phi,\zeta}f\|_{L^r(\R^n)}\le C\|N_{\Phi',\zeta'}f\|_{L^r(\R^n)}.\]
As a consequence, we can pick one and we can suppress the $\Phi,\zeta-$ dependence.\\

Lastly, recall that for any $1<r<\infty$ it holds \[\mathcal{H}^r(\R^n)=L^r(\R^n)\] with the equivalence of the norms, i.e. for any $1<r<\infty$ we can find constants $c_r>0$ and $C_r>0$ s.t. it holds
\begin{equation}c_r\|f\|_{L^r(\R^n)}\le\|f\|_{\mathcal{H}^r(\R^n)}\le C_r\|f\|_{L^r(\R^n)}.\label{eq:2.4.HS.equiv}\end{equation}
\end{de}

Secondly, we introduce the notion of a homogenous Lipschitz-space as can be found in \cite[Section 6.3]{G14}.

\begin{de}[Homogenous Lipschitz-spaces $\dot{\Lambda}_\beta$, $0<\beta\le1$] Define the semi-norm
\[\|f\|_{\dot{\Lambda}_\beta(\R^n)}:=\sup\limits_{\substack{x,y\in \R^n\\ x\not=y}}\frac{|f(x)-f(y)|}{|x-y|^\beta}.\]
Then we call
\[\dot{\Lambda}_\beta(\R^n):=\{f\in C(\R^n)|\|f\|_{\dot{\Lambda}_\beta(\R^n)}<\infty\}\]
\end{de}
 homogenous Lipschitz-space.\footnote{Note, that $\|\cd\|_{\dot{\Lambda}_\beta(\R^n)}$ is only a semi-norm on $\dot{\Lambda}_\beta(\R^n).$ However, $\dot{\Lambda}_\beta(\R^n)$ can be turned, as usual, into a normed space by $\dot{\Lambda}_\beta(\R^n)/\{constants\}.$}\\

\begin{re}1. Both types of spaces can be generalised to the case if the considered maps are vector-valued, as usual, we write $\mathcal{H}^r(\R^m,\R^n)$ and $\dot{\Lambda}_{\beta}(\R^m,\R^n).$ \\

2. For more information on Lipschitz spaces, for instance, a more general definition if $\be>1,$ or for the notion of inhomogeneous Lipschitz spaces see \cite[Section 6.3]{G14}.
\end{re}

Next, the duality result is stated, followed by a discussion of the literature, and where to find the proof.
\begin{lem} [$\mathcal{H}^r-\dot{\Lambda}_\beta$-Duality]
Assume $r<1,$ and $g\in \mathcal{H}^r(\R^n),$ and $f\in\dot{\Lambda}_\beta(\R^n)$ with $\beta:=n(\frac{1}{r}-1).$ Then there exists a constant $c(r,\Phi,\zeta)>0$ s.t.
\begin{equation}\int\limits_{\R^n}{fg\;dx}\le c(r,\Phi,\zeta)\|f\|_{\dot{\Lambda}_\beta(\R^n)}\|g\|_{\mathcal{H}^r(\R^n)}.
\label{eq:Ch1.4.103}\end{equation}
\end{lem}
\begin{re}
It is tough to find this statement in the literature.\footnote{Alternatively to the source we follow see \cite{JSW84}.} That is why we outline that the above is indeed true. We will follow Uchiyama's presentation, which he is using in his book \cite{U01} to stick to a single notation. \\

1. It has been shown that every function $f\in \mathcal{H}^r(\R^n)$ for $r\le1$ posses an atomic decomposition and hence $\mathcal{H}^r(\R^n)=\mathcal{H}_{at}^r(\R^n).$
Moreover, the norms are equivalent, too, i.e. there exists constant $c=c(r,\Phi,\zeta)>0$ and $C=C(r,\Phi,\zeta)>0$ s.t.
\[c(r,\Phi,\zeta)\|f\|_{\mathcal{H}_{at}^r(\R^n)}\le\|f\|_{\mathcal{H}^r(\R^n)} \le C(r,\Phi,\zeta)\|f\|_{\mathcal{H}_{at}^r(\R^n)},\]
see \cite[Eq.(0.8)-(0.9), p.7]{U01} for the statement and references to the original literature.\\

2. On the other hand Coifman and Weiss obtain in \cite[Theorem B, p.593]{Co77} the duality between $\mathcal{H}^r$ and $\dot{\Lambda}_\beta.$ 
It is important to realise that Coifman and Weiss use an equivalent notion to define these spaces by using a measure 
\[|f(x)-f(y)|\le \mu(B)^\alpha,\]
where $B$ is a ball containing $x,y$, see \cite[Eq.(2.2), p.591]{Co77}. By doing so obviously $\mu(B)\sim |x-y|^n$ and we get the relation $\beta=n\alpha.$\\

3. In particular, in \cite[Eq (3.19), p.634]{Co77} they obtain the duality inequality
\[\int\limits_{\R^n}{fg\;dx}\le \|f\|_{\dot{\Lambda}_\beta(\R^n)}\|g\|_{\mathcal{H}_{at}^r(\R^n)}.\]
This is indeed the case, however, they use a slightly different notation compared to Uchiyama. Coifman and Weiss use $|g|_{r,1}^{1/r}$ while Uchiyama has already absorbed the the exponent $1/p$ in the definition of $\|g\|_{\mathcal{H}_{at}^r(\R^n)},$ yielding the same result. Moreover, with the equivalence of the norms in point 1 we get 
\[\int\limits_{\R^n}{fg\;dx}\le C(r,\d,\phi)\|f\|_{\dot{\Lambda}_\beta(\R^n)}\|g\|_{\mathcal{H}^r(\R^n)},\]
which agrees with the conjecture.
\end{re}
We end this section by stating the famous div-curl lemma first shown in \cite{CLMS93}.
\begin{lem}[div-curl lemma]\label{lem:2.HR.1}
Let $\frac{n}{n+1}<s,t<\infty$ s.t. $\frac{1}{s}+\frac{1}{t}<1+\frac{1}{n},$ and $s, t$ are such that at least one of them is strictly larger than 1, and assume $F\in \mathcal{H}^s(\R^n,\R^n)$  and $G\in \mathcal{H}^t(\R^n,\R^n)$ satisfy \[\mb{curl}\, F=0 \;\mb{in}\; D'(\R^n,\R^\sigma)\;\mb{and}\; \div G=0 \;\mb{in}\; D'(\R^n),\] where $\sigma:=\mb{dim}\!\mb{curl}\!(\mathcal{H}^s(\R^n,\R^n)).$ Then $F\cd G\in \mathcal{H}^r(\R^n)$ for $\frac{1}{r}=\frac{1}{s}+\frac{1}{t}$ with
\begin{equation}\|F\cd G\|_{\mathcal{H}^r(\R^n)}\le C\|F\|_{\mathcal{H}^s(\R^n,\R^n)}\|G\|_{\mathcal{H}^t(\R^n,\R^n)},\end{equation}
for some positive constant $C>0.$
\end{lem}
\begin{re} Recall, that it holds
\[\div G=0 \;\mb{in}\; D'(\R^n) \;\mb{iff} \;\int\limits_{\R^n}{G\cd\grad h\;dx}=0 \;\mb{for any}\; h\in C_c^\infty(\R^n)\]
for any dimension $n\in \N\sm\{0\}.$\\
The distributional version of 'curl' is more complicated. For the relevant dimension $n=2$ it holds
\[\mb{curl}\, F=0 \;\mb{in}\; D'(\R^2) \;\mb{iff} \;\int\limits_{\R^2}{(F_2h,_1-F_1h,_2)\;dx}=0 \;\mb{for any}\; h\in C_c^\infty(\R^2). \]
A similar version can be obtained for $n=3.$ In higher dimensions one might want to use differential forms to define a distributional version of the 'curl'. The details are omitted for simplicity.
\end{re}

\bibliography{LiteraturePhDMD}
\bibliographystyle{plain}
\vspace{0.5cm}
\textsc{Acknowledgements:} The author is appreciative to the Department of Mathematics at the University of Surrey and was funded by the Engineering \& Physical Sciences Research Council (EPRSC). Thank you to Jonathan J. Bevan for suggesting crucial parts of the proofs of Theorem 1.2 and 1.4 and also to Bin Cheng for comments.
\end{document}